\documentclass[11pt]{amsart}

 \usepackage{mathpazo}
\usepackage[scaled=0.90]{helvet} 
\usepackage{courier} 
\normalfont
\usepackage[T1]{fontenc}
\usepackage[sort]{cite} 

\numberwithin{equation}{section} \oddsidemargin=-.0cm
\usepackage{textcomp}
\usepackage[ruled]{algorithm2e}
\usepackage[usenames, dvipsnames]{color}

\usepackage[colorlinks=true, pdfstartview=FitV, linkcolor=blue,
            citecolor=blue, urlcolor=blue]{hyperref}

\numberwithin{equation}{section}
\usepackage{framed,dsfont}
\usepackage{amsfonts, amssymb, amsmath, graphics, subfigure,mathtools,amsbsy}
\usepackage{booktabs}

\graphicspath{{./Figures/}}

\newtheorem{thm}{Theorem}[section]
\newtheorem{lem}{Lemma}[section]
\newtheorem{defi}{Definition}[section]

\newtheorem{rem}{Remark}[section]
\newtheorem{cor}{Corollary}[section]
\def\bt{\begin{thm}}
\def\et{\end{thm}}
\def\bl{\begin{lem}}
\def\el{\end{lem}}
\def\bd{\begin{defi}}
\def\ed{\end{defi}}
\def\bc{\begin{cor}}
\def\ec{\end{cor}}
\def\bp{\begin{proof}}
\def\ep{\end{proof}}
\def\br{\begin{rem}}
\def\er{\end{rem}}
\def\bi{\begin{itemize}}
\def\ei{\end{itemize}}

\def\be{\begin{equation}}
\def\ee{\end{equation}}
\def\bes{\begin{equation*}}
\def\ees{\end{equation*}}
\def\bea{\begin{equation} \begin{aligned}}
\def\eea{\end{aligned} \end{equation}}
\def\beas{\begin{equation*} \begin{aligned}}
\def\eeas{\end{aligned} \end{equation*}}

\def\Hlf{H_{\text{LF}}}

\def\tr{\mathrm{tr}}
\def\Lip{{\rm Lip}}

\def\d{\mathrm{d}}
\def\cH{\mathcal H}
\def\Forall{\text{ } \forall \:}
\def\R{\mathbb{R}}
\def\Z{\mathbb{Z}}
\def\N{\mathbb{N}}

\newcommand{\mk}{\color{black}}

\title[Galerkin approximations of optimal control of nonlinear DDEs]{Galerkin approximations for the optimal control of nonlinear delay differential equations}

\author[Micha\"el D. Chekroun]{Micka\"el D. Chekroun}
\address{Department of Atmospheric \& Oceanic Sciences, University of California, Los Angeles, CA 90095-1565, USA}
\email{mchekroun@atmos.ucla.edu}

\author[Axel Kr\"oner]{Axel Kr\"oner}
\address{Institut f\"ur Mathematik, Humboldt-Universit\"at Berlin, 10099 Berlin, Germany; INRIA and CMAP, \'Ecole Polytechnique, CNRS, Universit\'e Paris Saclay, 91128 Palaiseau, France}
\email{axel.kroener@inria.fr}

\author[Honghu Liu]{Honghu Liu}
\address{Department of Mathematics, Virginia Polytechnic Institute and State University, Blacksburg, Virginia 24061, USA}
\email{hhliu@vt.edu}


\begin{document}

\maketitle

\begin{abstract}
Optimal control problems of nonlinear delay differential equations (DDEs) are considered for which we propose a general 
Galerkin approximation scheme built from Koornwinder polynomials. Error estimates for the resulting Galerkin-Koornwinder approximations to the optimal control and the value function, are derived for a broad class of cost functionals and nonlinear DDEs. The approach is illustrated on a delayed logistic equation set not far away from its Hopf bifurcation point in the parameter space. In this case, we show that low-dimensional controls for a standard quadratic cost functional can be efficiently computed from Galerkin-Koornwinder approximations to reduce at a nearly optimal cost the oscillation amplitude displayed by the DDE's solution. Optimal controls computed from the Pontryagin's maximum principle (PMP) and the Hamilton-Jacobi-Bellman equation (HJB) associated with the corresponding ODE systems, are shown to provide numerical solutions in good agreement. It is finally argued that the value function computed from the  corresponding reduced HJB equation provides a good approximation of that obtained from the full HJB equation.
\end{abstract}

\tableofcontents

\section{Introduction}

Delay differential equations (DDEs) are widely used in many fields such as biosciences \cite{Diekmann_al95, Kuang93,Smith11,MacDonald89},
climate dynamics \cite{CGN17,GCStep15, Sieber2014,Suarez_al88,TSCJ94}, chemistry and engineering \cite{Hale_Lunel93, LS10,  Michiels_al07, Stepan89,Kolmanovskii_al1986}.  The inclusion of time-lag terms are aimed to account for delayed responses of the modeled systems to either internal or external factors. Examples of such factors include incubation period of infectious diseases \cite{Kuang93}, wave propagation \cite{Chek_Glatt-Holtz,Suarez_al88}, or time lags arising in engineering \cite{Kolmanovskii_al1986} to name a few.

In contrast to ordinary differential equations (ODEs), the phase space associated even with a scalar DDE is infinite-dimensional, due to the presence of delay terms implying thus the knowledge of a function as an initial datum, namely the initial history to solve the Cauchy problem \cite{curtain1995,Hale_Lunel93}; cf.~Section~\ref{Sect_theory}. The infinite-dimensional nature of the state space renders the related optimal control problems more challenging to solve compared to the ODE case. It is thus natural to seek for low-dimensional approximations to alleviate the inherent computational burden. The need of such approximations is particularly relevant when (nearly) optimal controls in feedback form are sought, to avoid solving the infinite-dimensional Hamilton-Jacobi-Bellman (HJB) equation  associated with the optimal control problem of the original DDE.

For optimal control of {\it linear} DDEs, averaging methods as well as spline approximations are often used for the design of (sub)optimal controls in either open-loop form or feedback form; see e.g.~\cite{Banks_al78,Banks_al1984,Kappel_al1987,Kunisch82} and references therein. The usage of spline approximations in the case of open-loop control for {\it nonlinear} DDEs has also been considered \cite{Lamm84}, but not in a systematic way. Furthermore, 
due to the locality of the underlying basis functions, the use of e.g.~spline approximations for the design of feedback controls leads, especially in the nonlinear case, to surrogate HJB equations of too high dimension to be of  practical interest.

In this article, we bring together a recent approach dealing with the Galerkin approximations for  optimal control problems of general nonlinear evolution equations in a Hilbert space \cite{CKL17}, with  techniques for the finite-dimensional and analytic approximations of nonlinear DDEs based on Legendre-type polynomials, namely the Koornwinder polynomials \cite{CGLW15}.  Within the resulting framework, we adapt ideas from \cite[Sect.~2.6]{CKL17}  to derive---for a broad class of cost functionals and nonlinear DDEs---error estimates for the approximations of the value function and optimal control, associated with the Galerkin-Koornwinder (GK)  systems; see Theorem~\ref{Thm_PM_val} and Corollary~\ref{Lem_controller_est}. These error estimates are formulated in terms of residual energy, namely the energy contained in the controlled DDE solutions as projected onto the orthogonal complement of a given Galerkin subspace.

Our approach is then applied to a Wright equation\footnote{Analogue to a logistic equation with delay; see Sect.~\ref{Sect_Wright_PMP}.} with the purpose of  reducing (optimally) the amplitude of the oscillations displayed by the DDE's solution subject to a quadratic cost; see Sect.~\ref{Sect_DDE}. For this model, the oscillations emerge through a supercritical Hopf bifurcation as the delay parameter $\tau$ crosses a critical value, $\tau_c$, from below. We show, for $\tau$ above and close to  $\tau_c$, that a suboptimal controller at a nearly optimal cost can be synthesized from a 2D projection of the 6D GK system; see Sect.~\ref{Sec_num_results}. The projection is made onto the space spanned by the first two eigenvectors of the linear part of the 6D GK system; the sixth dimension constituting for this example the minimal dimension---using Koornwinder polynomials---to resolve accurately the linear contribution to the oscillation frequency of the DDE's solution; see \eqref{Eq_osc_part} in Sect.~\ref{Sec_num_results}.

Using the resulting 2D ODE system, the syntheses of (sub)optimal controls obtained either from application of the Pontryagin's maximum principle (PMP) or by solving the associated HJB equation, are shown to provide nearly identical numerical solutions. Given the good control skills obtained from the controller synthesized from the reduced HJB equation, one can reasonably infer that the corresponding ``reduced'' value function  provides thus a good approximation of the ``full'' value function associated with the DDE optimal control problem; see Sect.~\ref{Sec_HJB}.

The article is organized as follows. In Sect.~\ref{Sect_theory}, we first introduce the class of nonlinear DDEs considered in this article and recall in Sect.~\ref{Sec_DDEreformulation} how to recast such DDEs into infinite--dimensional evolution equations in a Hilbert space. We summarize then in Sect.~\ref{Sect_GK_approx} the main tools from \cite{CGLW15}  for the analytic determination of GK systems approximating DDEs. In Sect.~\ref{Sect_DDE_control},  we  derive error estimates for the approximations of the value function and optimal control, obtained from the GK  systems. Our approach is applied to the Wright equation in Sect.~\ref{Sect_DDE}. Section \ref{Sect_Wright_PMP} introduces the 
optimal control problem associated with this equation and its  abstract functional formulation following Sect.~\ref{Sec_DDEreformulation}. 
 The explicit GK approximations of the corresponding optimal control problem are derived in Sect.~\ref{Sect_Wright_GK_control}. Numerical results based on PMP from a projected 2D GK system are then presented in Sect.~\ref{Sec_num_results}. 
In Sect.~\ref{Sec_HJB}, we show that by solving numerically the associated reduced HJB equation, a (sub)optimal control in a feedback form can be synthesized at a nearly optimal cost.  Finally,  directions to derive reduced systems of even lower dimension from GK approximations are outlined in  Sect.~\ref{Sect_conclusion}. 
 

\section{Optimal control of DDEs: Galerkin-Koornwinder approximations} \label{Sect_theory}

In this article, we are concerned with optimal control problems associated with nonlinear DDEs of the form
\be \label{Eq_DDE}
\frac{\d m(t)}{\d t} = a m(t) + b m(t-\tau) + c \int_{t-\tau}^t m(s)\d s + F\Big(m(t), m(t-\tau),\int_{t-\tau}^t m(s) \d s\Big),
\ee
where $a$, $b$ and $c$ are real numbers, $\tau> 0$ is the delay parameter, and $F$
is a nonlinear function. We refer to Sect.~\ref{Sect_DDE_control} for assumption on $F$. To simplify the presentation, we restrict ourselves to this class of scalar DDEs, but nonlinear systems of DDEs involving several delays are allowed by the approximation framework of \cite{CGLW15} adopted in this article.

Note that even for scalar DDEs such as given by Eq.~\eqref{Eq_DDE}, the associated state space is infinite-dimensional. This is due to the presence of time-delay terms, which require providing initial data over an interval $ [-\tau,0],$ where $\tau >0$ is the delay. 
It is often desirable, though, to have low-dimensional ODE systems that capture qualitative features,  as well as approximating certain quantitative aspects of the DDE dynamics. 

The approach adopted here consists of approximating the infinite-dimensional equation \eqref{Eq_DDE} by a finite-dimensional ODE system built from Koornwinder polynomials {\mk following} \cite{CGLW15}, and then solve reduced-order optimal control problems aimed at approximating a given optimal control problem associated with Eq.~\eqref{Eq_DDE}. 

To justify this approach, Eq.~\eqref{Eq_DDE} needs first to be recast into an evolution equation, in order to apply in a second step, recent results dealing with the {\mk Galerkin approximations to optimal control problems governed 
by general nonlinear evolution equations in Hilbert spaces, such as  presented in \cite{CKL17}.}

As a cornerstone, rigorous Galerkin approximations to Eq.~\eqref{Eq_DDE} are crucial, and are recalled hereafter.  We first recall how a DDE such as Eq.~\eqref{Eq_DDE} can be recast into an infinite-dimensional evolution equation in a Hilbert space.

\subsection{Recasting a DDE into an infinite--dimensional evolution equation}\label{Sec_DDEreformulation}
The reformulation of a system of DDEs into an infinite-dimensional ODE is classical. For this purpose, two types of function spaces are typically used  as state space:  {\mk the space of continuous functions $C([-\tau, 0]; \mathbb{R}^d)$, see \cite{Hale_Lunel93}, and the Hilbert space $L^2([-\tau, 0); \mathbb{R}^d)\times \mathbb{R}^d$, see \cite{curtain1995}.  The Banach space setting of continuous functions has been extensively used in the study of bifurcations arising from DDEs, see e.g.~\cite{Casal_al80,Chow_al77,Das_al02,Faria_al95,Kazarinoff_al78,Nayfeh08,wischert1994delay}, while the Hilbert space setting is typically adopted for the approximation of DDEs or their optimal control; see e.g.~\cite{Banks_al78,banks1979spline,Banks_al84,CGLW15,kappel1978autonomous,Kappel86,Kappel_al87,Kunisch82,Federico_al10,Ito_Teglas86}.}

Being concerned with the optimal control of scalar DDEs of the form \eqref{Eq_DDE}, we adopt here the Hilbert space setting and consider in that respect our state space to be
\be \label{H_space}
\mathcal{H} := L^2([-\tau,0); \mathbb{R}) \times  \mathbb{R},
\ee 
endowed with the {\mk following inner product,} defined for any  $(f_1, \gamma_1)$ and $(f_2, \gamma_2)$  in $\mathcal{H}$ by:
\be \label{H_inner}
 \langle (f_1, \gamma_1), (f_2, \gamma_2) \rangle_{\mathcal{H}} := \frac{1}{\tau} \int_{-\tau}^0 f_1(\theta) f_2(\theta)  \d \theta  + \gamma_1 \gamma_2.
\ee

We will also make use sometimes of the following subspace of $\mathcal{H}$:
\be
\mathcal{W} := H^1([-\tau,0); \mathbb{R}) \times  \mathbb{R},
\ee
where $H^1([-\tau,0); \mathbb{R})$ denotes the standard Sobolev subspace of $L^2([-\tau,0); \mathbb{R})$.

Let us denote by $m_t$ the time evolution of the history segments of a solution to Eq.~\eqref{Eq_DDE}, namely 
\be \label{shift}
 m_t(\theta):=m(t+\theta), \qquad t \ge 0, \qquad \theta \in [-\tau, 0].
\ee
Now, by introducing the new variable
\be
y(t) := (m_t,m(t))=(m_t, m_t(0)), \; t\geq 0,
\ee
Eq.~\eqref{Eq_DDE} can be rewritten as the following  abstract ODE in the Hilbert space $\cH$: 
\be \label{eq:DDE_abs0}
\frac{\d y}{\d t} = \mathcal{A} y + \mathcal{F}(y),
\ee
where the linear operator $\mathcal{A} \colon D(\mathcal{A}) \subset \mathcal{W} \rightarrow \mathcal{H}$ is defined as
\bea \label{Def_A}
\lbrack \mathcal{A} \Psi \rbrack (\theta) & := \begin{cases}
{\displaystyle \frac{\d^+ \Psi^D}{\d \theta}}, &  \theta \in[-\tau, 0),  \vspace{0.4em}\\ 
{\displaystyle a \Psi^S + b\Psi^D(-\tau) + c \int_{-\tau}^0 \Psi^D(s)\d s}, & \theta = 0,
\end{cases} 
\eea
with the domain $D(\mathcal{A})$ of $\mathcal{A}$ given by (cf.~e.g.~\cite[Prop.~2.6]{Kappel86})
\be \label{D_of_A}
D(\mathcal{A}) = \Big \{(\Psi^D, \Psi^S) {\mk \in \cH}: \Psi^D \in H^1([-\tau, 0); \mathbb{R}), \lim_{\theta \rightarrow 0^-} \Psi^D(\theta) = \Psi^S 
\Big \}.
\ee
The nonlinearity $\mathcal{F} \colon \mathcal{H} \rightarrow \mathcal{H}$ is here {\mk given, for all $\Psi = (\Psi^D, \Psi^S)$ in  $\mathcal{H}$,}
 by 
\bea \label{Def_F}
[\mathcal{F} (\Psi) ](\theta) & := \begin{cases}
0, &  \theta \in[-\tau, 0),   \vspace{0.4em}\\ 
F \Big(\Psi^S, \Psi^D(-\tau), \int_{-\tau}^0 \Psi^D(s) \d s \Big), & \theta = 0, 
\end{cases}  \eea

With $D(\mathcal{A})$ such as given in \eqref{D_of_A}, the operator $\mathcal{A}$ generates a linear $C_0$-semigroup on $\mathcal{H}$   so that the Cauchy problem associated with the linear  equation  $\dot{y}=\mathcal{A} y$ is well-posed in the Hadamard's sense; see e.g~\cite[Thm.~2.4.6]{curtain1995}. The well-posedness problem for the nonlinear equation depends obviously on the nonlinear term $\mathcal{F}$ and we refer to \cite{CGLW15} {\mk and references therein for a discussion of this problem}; see also \cite{Webb76}. 

For later usage, we recall that a {\it mild solution} to Eq.~\eqref{eq:DDE_abs0} over $[0, T]$ with initial datum $y(0) = x$ in $\cH$, {\mk is an element $y$ in $C([0,T],\cH)$ that satisfies the integral equation}
\be \label{Def_mild_soln}
y(t)=T(t)x + \int_0^t T(t-s) \mathcal{F}(y(s)) \d s, \qquad \forall\; t  \in [0,T],
\ee
 where {\mk $(T(t))_{t\geq 0}$} denotes the $C_0$-semigroup generated by the operator 
 $\mathcal{A}\colon D(\mathcal{A}) \rightarrow \mathcal{H}$. {\mk This notion of mild solutions extend naturally when a control term $\mathfrak{C}(u(t))$ is added to the RHS of Eq.~\eqref{eq:DDE_abs0}; in such a case the definition of a mild solution is amended by the presence of the integral term  $\int_0^t T(t-s) \mathfrak{C}(u(s)) \d s$ to the RHS of Eq.~\eqref{Def_mild_soln}.}

\subsection{Galerkin-Koornwinder approximation of DDEs} \label{Sect_GK_approx}

{\mk Once a DDE is reframed into an infinite-dimensional ODE in $\cH$, the approximation problem by finite-dimensional ODEs can be addressed in various ways.  In that respect,} different basis functions have been proposed to decompose  
the state space $\cH$;  these include, among others, step functions \cite{Banks_al78,kappel1978autonomous}, splines \cite{banks1979spline,Banks_al84}, and orthogonal  polynomial functions, such as Legendre polynomials \cite{Kappel86,Ito_Teglas86}. Compared with step functions or splines, the use of orthogonal polynomials leads typically to ODE approximations with lower dimensions, for a given precision \cite{banks1979spline,Ito_Teglas86}. On the other hand, classical polynomial basis functions  do not live in the domain of the linear operator underlying the DDE, which leads to technical complications in establishing convergence results \cite{Kappel86,Ito_Teglas86}; see \cite[Remark 2.1-(iii)]{CGLW15}.

{\mk
One of the  main contributions of \cite{CGLW15} consisted in identifying that Koornwinder polynomials \cite{Koo84} lie in the domain $D(\mathcal{A})$ of linear operators such as $\mathcal{A}$ given in \eqref{Def_A}, allowing in turn for adopting a classical Trotter-Kato (TK) approximation approach from the $C_0$-semigroup theory \cite{Goldstein85,Pazy83} to deal with the ODE approximation of DDEs such as given by Eq.~\eqref{Eq_DDE}. 
The TK approximation approach can be viewed as the functional analysis operator version of the Lax equivalence principle.\footnote{i.e., if ``consistency'' and ``stability'' are satisfied, then ``convergence'' holds, and reciprocally.}  The work \cite{CGLW15} allows thus for positioning the approximation problem of DDEs within a well defined territory. In particular, as pointed out in \cite{CKL17} and discussed in Sect.~\ref{Sect_DDE_control} below, the optimal control of DDEs benefits from the TK approximation approach. 

In this section, we focus on another important feature pointed out in \cite{CGLW15} for applications, namely, Galerkin approximations of DDEs built from Koornwinder polynomials can be efficiently computed via simple analytic formulas; see \cite[Sections 5-6 and Appendix C]{CGLW15}.
We recall below the main elements to do so referring to \cite{CGLW15} for more details.} 

{\mk First, let us recall that Koornwinder polynomials $K_n$ are obtained from 
Legendre polynomials $L_n$ according to the relation}
\be \label{eq:Pn}
K_n(s) := -(1+s)\frac{\d}{\d s} L_n(s) +( n^2 + n + 1) L_n(s), \; s \in [-1, 1], \; n \in \mathbb{N};
\ee
{\mk see \cite[Eq.~(2.1)]{Koo84}.}

Koornwinder polynomials are  known to form an orthogonal set for the following weighted inner product with a point-mass on $[-1,1]$,
\be\label{Eq_Koorn_dx}
\mu(\d x)= \frac{1}{2} \d x  + \delta_{1}, 
\ee 
where $\delta_1$ denotes the Dirac point-mass at the right endpoint $x=1$; see  \cite{Koo84}.  In other words,
\bea
\int_{-1}^{1} K_n(s) K_p(s) \d \mu (s)& = \frac{1}{2} \int_{-1}^{1} K_n(s) K_p(s) \d s + K_n(1) K_p(1)\\
&=0, \, \mbox{ if } p\neq n.
\eea

It is also worthwhile noting that  the sequence given by
\be \label{eq:Pn_prod}
\{\mathcal{K}_n := (K_n, K_n(1)) : n \in \mathbb{N}\}
\ee 
forms an orthogonal basis of the  product space 
\be \label{eq:E}
\mathcal{E} := L^2([-1,1); \mathbb{R}) \times  \mathbb{R},
\ee 
where $\mathcal{E}$ is endowed with the following inner product:
\be \label{eq:inner_E}
\langle (f, a), (g, b) \rangle_{\mathcal{E}}  = \frac{1}{2} \int_{-1}^1 f(s)g(s) \d s  + ab, \quad (f,a), (g, b) \in \mathcal{E}.
\ee
The norm induced by this inner product will be denoted hereafter by $\|\cdot\|_{\mathcal{E}}$.

From the original Koornwinder basis given on the interval $[-1, 1]$, orthogonal polynomials on the interval $[-\tau, 0]$ for the inner product \eqref{H_inner} can now be easily obtained by using a simple linear transformation $\mathcal{T}$ defined by:
\be \label{eq:linear_transf}
\mathcal{T} \colon [-\tau, 0] \rightarrow [-1, 1], \qquad \theta \mapsto 1 + \frac{2 \theta }{\tau}. 
\ee
Indeed, for $K_n$ given by \eqref{eq:Pn}, let us define the {\mk rescaled} polynomial $K_n^\tau$ by 
\bea \label{eq:Pn_tilde}
K^\tau_n\colon  [-\tau, 0] & \rightarrow \mathbb{R}, \\
\theta & \mapsto  K_n \Bigl( 1 + \frac{2 \theta }{\tau} \Bigr), \qquad n \in \mathbb{N}.
\eea

As shown in \cite{CGLW15}, the sequence  
\be \label{eq:Pn_tilde_prod}
\{\mathcal{K}_n^\tau := (K_n^\tau, K_n^\tau(0)) : n \in \mathbb{N}\}
\ee
forms an orthogonal basis for the space $\mathcal{H} = L^2([-\tau,0); \mathbb{R}) \times  \mathbb{R}$ endowed with the inner product $\langle \cdot, \cdot \rangle_{\mathcal{H}}$ given in \eqref{H_inner}. Note that since $K_n(1)=1$ \cite[Prop.~3.1]{CGLW15}, we have 
\be\label{Eq_normalization}
K_n^\tau(0)  = 1.
\ee

As shown in \cite[Sect.~5 \& Appendix C]{CGLW15}, (rescaled) Koornwinder polynomials allow for analytical Galerkin approximations of general nonlinear systems of DDEs. 
In the case of a nonlinear scalar DDE such as Eq.~\eqref{Eq_DDE}, the $N^{\mathrm{th}}$ Galerkin-Koornwinder (GK) approximation, $y_N$, is obtained as
\be \label{x_t expand}
y_N(t) = \sum_{j=0}^{N-1} \xi^{N}_j(t) \mathcal{K}^\tau_j, \qquad t \ge 0,
\ee
where the $\xi^N_j(t)$ solve the $N$-dimensional ODE system 
\begin{equation} \label{Galerkin_AnalForm}
\begin{aligned}
\frac{\d \xi^N_j}{\d t} & = \frac{1}{\|\mathcal{K}_j\|_{\mathcal{E}}^2 } \sum_{n=0}^{N-1} \Big( a + b K_n(-1) + c \tau (2 \delta_{n,0} - 1) \\
& \hspace{8em} + \frac{2}{\tau}\sum_{k=0}^{n-1} a_{n,k} \left( \delta_{j,k} \|\mathcal{K}_j\|^2_{\mathcal{E}} - 1 \right) \Big) \xi^N_n(t) \\
& \hspace{1em} + \frac{1}{\|\mathcal{K}_j\|_{\mathcal{E}}^2} F \left( \sum_{n=0}^{N-1} \xi^N_n(t),  \sum_{n=0}^{N-1} \xi^N_n(t) K_n(-1),\tau \xi^N_0(t) - \tau \sum_{n=1}^{N-1} \xi^N_n(t) \right),\\
& \;\; 0\leq j\leq N-1.
\end{aligned}
\end{equation}
{\mk Here the Kronecker symbol $\delta_{j,k}$ has been used,} and  the coefficients $a_{n,k}$ are obtained by solving a triangular linear system in which the right hand side has 
explicit coefficients depending on $n$; see \cite[Prop.~5.1]{CGLW15}.

In practice, an approximation $m_N(t)$ of $m(t)$ solving Eq.~\eqref{Eq_DDE} is obtained   
 as the state part (at time $t$) of $y_N$ given by  \eqref{x_t expand} which, thanks to the normalization property $K_n^\tau(0)  = 1$ given in ~\eqref{Eq_normalization}, reduces to
\be \label{eq_mN}
m_N(t) = \sum_{j = 0}^{N-1} \xi^N_j(t).
\ee


For later usage, we rewrite the above GK system into the following compact form:
\be \label{Galerkin_cptForm}
\boxed{\frac{\d \boldsymbol{\xi}_N}{\d t} = M \boldsymbol{\xi}_N + G (\boldsymbol{\xi}_N),}
\ee
where $M\boldsymbol{\xi}_N$ denotes the linear part of Eq.~\eqref{Galerkin_AnalForm}, and $G(\boldsymbol{\xi}_N)$ the nonlinear part. Namely, $M$ is the $N\times N$ matrix whose elements are given by
\begin{equation} \label{eq:A}
\boxed{
\begin{aligned}
(M)_{j,n} & = \frac{1}{\|\mathcal{K}_j\|_{\mathcal{E}}^2 }\Big(a + b K_n(-1) + c \tau (2 \delta_{n,0} - 1) \\
& \hspace{8em}+ \frac{2}{\tau}\sum_{k=0}^{n-1} a_{n,k} \left( \delta_{j,k} \|\mathcal{K}_j\|^2_{\mathcal{E}} - 1 \right ) \Big),
\end{aligned}
}
\end{equation}
where $j, n = 0, \cdots, N-1$, and the nonlinear vector field $G \colon \mathbb{R}^N \rightarrow \mathbb{R}^N$, is given component-wisely
by
\be \label{eq:G}
\boxed{G_j(\boldsymbol{\xi}_N) = \frac{1}{\|\mathcal{K}_{j}\|_{\mathcal{E}}^2} F \left( \sum_{n=0}^{N-1} \xi^N_n(t),  \sum_{n=0}^{N-1} \xi^N_n(t) K_n(-1),\tau \xi^N_0(t) - \tau \sum_{n=1}^{N-1} \xi^N_n(t) \right),}
\ee
for each  $0\leq j\leq N-1$.

 For rigorous convergence results of GK systems \eqref{Galerkin_cptForm} for Eq.~\eqref{Eq_DDE}
 when $F$ does not depend on the delay term $m(t-\tau)$,  we refer to \cite[Sect.~4]{CGLW15}.  {\mk When $F$ does depend on $m(t-\tau)$,  convergence is also believed to hold, as supported by the striking numerical results shown in  \cite[Sect.~6]{CGLW15}.}

\subsection{Galerkin-Koornwinder approximation of nonlinear optimal control problems associated with DDEs} \label{Sect_DDE_control}

\subsubsection{Preliminaries} 

In this section, we outline how the recent results of \cite{CKL17} concerned with the 
Galerkin approximations to optimal control problems governed by nonlinear evolution equations in Hilbert spaces, apply to the context of nonlinear DDEs when these approximations are built from the Koornwinder polynomials. We provide thus here further elements regarding the research program outlined in \cite[Sect.~4]{CKL17}.

We consider here, given a finite horizon $T>0$, the following controlled version of Eq.~\eqref{eq:DDE_abs0},
\bea \label{eq:DDE_abs_c}
\frac{\d y}{\d s}& = \mathcal{A} y + \mathcal{F}(y)+\mathfrak{C} (u(s)),
\quad s \in (t, T], \; u\in \mathcal{U}_{ad}[t,T],\\
y(t) &=x \, \in \cH.
\eea
The (possibly nonlinear) mapping $\mathfrak{C}: V \rightarrow \cH$ is assumed to be such that $\mathfrak{C}(0)=0$, and the control $u(s)$ {\mk lies} in a bounded subset $U$ of a separable Hilbert space $V$ possibly different from $\cH$.
Assumptions about the set, $\mathcal{U}_{ad}[t,T]$, of admissible controls, $u$,  is made precise below; see \eqref{Eq_U_bounded}. 

Here our Galerkin subspaces are spanned by the rescaled Koornwinder polynomials, namely
\be
\cH_N = \mathrm{span}\{\mathcal{K}^\tau_0, \cdots, \mathcal{K}^\tau_{N-1}\}, \quad N = 1,2,\cdots,
\ee
where $\mathcal{K}^\tau_N$ is defined in \eqref{eq:Pn_tilde_prod}. {\mk Denoting} by $\Pi_N$ the projector of $\cH$ onto $\cH_N$, the corresponding GK approximation of Eq.~\eqref{eq:DDE_abs_c} is then {\mk given by}
\bea \label{eq:DDE_Galerkin_abs}
\frac{\d y_N}{\d s}&= \mathcal{A}_N y_N + \Pi_N\mathcal{F}(y_N)+\Pi_N\mathfrak{C} (u(s)), \quad s \in (t, T], \; u\in \mathcal{U}_{ad}[t,T],\\
y_N(t)&=\Pi_N x \in \cH_N,
\eea
where $\mathcal{A}_N:=\Pi_N \mathcal{A} \Pi_N$.

Given a cost functional assessed along a solution
$y_{t,x}(\cdot;u)$ of Eq.~\eqref{eq:DDE_abs_c} driven by $u$
\be  \label{J_tx}
J_{t,x}(u) := \int_t^T [\mathcal{G}(y_{t,x}(s;u)) + \mathcal{E}(u(s)) ] \d s, \qquad t \in [0,T), \; u \in \mathcal{U}_{ad}[t,T],
\ee
{\mk (with $\mathcal{G}$ and $\mathcal{E}$ to be defined)} the focus of \cite{CKL17} was to identify simple checkable conditions that guarantee in---but not limited to---such a context, 
the convergence of the corresponding value functions, namely  
\be \label{value_est_goal}
\lim_{N \rightarrow \infty} \sup_{t \in [0, T]} |v_N(t,\Pi_N x) - v(t,x)| = 0,
\ee
where
\begin{subequations}\label{Eq_val_fcts_sec3}
\begin{align}
& v(t, x) := \inf_{u \in \mathcal{U}_{ad}[t,T]} J_{t,x}(u),\;  (t,x) \in [0,T) \times \cH \;  \text{ and } \;  v(T, x) = 0, \label{subEq1}\\
& v_N(t, x_N) := \inf_{u \in \mathcal{U}_{ad}[t,T]} J^N_{t,x_N}(u), \; (t,x_N) \in [0,T)\times \cH_N \;  \text{ and } \; v_N(T, x_N) = 0. \label{subEq2}
\end{align}
\end{subequations}
and $J^N_{t,x_N}(u)$ denotes the cost functional $J$ assessed along the Galerkin approximation  $y^N_{t,x_N}(\cdot;u)$ to \eqref{eq:DDE_abs_c} driven by $u$ and emanating at time $t$ from
\be
x_N:=\Pi_N x.
\ee

As shown in \cite{CKL17}, conditions ensuring {\mk the convergence \eqref{value_est_goal}} can be grouped {\mk into} three categories. The first set of conditions deal with the operator $\mathcal{A}$ and its Galerkin approximation $\mathcal{A}_N$. Essentially, it boils down to show that $\mathcal{A}$ generates a $C_0$-semigroup on $\cH$, and that the Trotter-Kato approximation conditions are satisfied; see Assumptions {\bf (A0)}--{\bf (A2)} in \cite[Section 2.1]{CKL17}. The latter conditions are, {\mk as mentioned earlier,} satisfied in the case of DDEs and GK approximations as shown by \cite[Lemmas 4.2 and 4.3]{CGLW15}. 

The second group of conditions identified in \cite{CKL17} is also not restrictive in the sense that only local Lipschitz  conditions\footnote{Note however that in the case of DDEs, nonlinearities $\mathcal{F}$ that include discrete delays may complicate the verification of a local Lipschitz property in $\cH$ as given by \eqref{H_space}. The case of  nonlinearities $\mathcal{F}$  depending on distributed delays and/or the current state can be handled in general more easily; see \cite[Sect.~4.3.2]{CGLW15}.} on $\mathcal{F}$, $\mathcal{G}$,  and $\mathfrak{C}$ are required as well as continuity of $\mathcal{E}$ and compactness of $U$ where $\mathcal{U}_{ad}[t,T]$ is taken to be
\be\label{U_ad_toto}
\mathcal{U}_{ad}[t,T]:=\{u\vert_{[t,T]} \; : \; u \in \mathcal{U}_{ad}\},
\ee
with 
\be\label{Eq_U_bounded}
\mathcal{U}_{ad}:=\{u\in L^q(0,T; V)\;:  \; u(s) \in U \textrm{ for a.e. } s \in [0,T]\},\; \;q\geq 1.
\ee
Also required to hold, are a priori bounds---uniform in $u$ in $\mathcal{U}_{ad}$---for the solutions of Eq.~\eqref{eq:DDE_abs_c} as well as of their Galerkin approximations. Depending on the specific form of Eq.~\eqref{eq:DDE_abs_c} such bounds can be derived for a broad class of DDEs; in that respect the proofs of \cite[Estimates (4.75)]{CGLW15} and \cite[Corollary 4.3]{CGLW15} can be easily adapted to the case of controlled DDEs. 

Finally, the last condition to ensure \eqref{value_est_goal}, requires that the residual energy of the solution $s\mapsto y(s;x,u)$ of 
\eqref{eq:DDE_abs_c} (with $y(0)=x$) driven by $u$, satisfies 
\be
\|\Pi^\perp_N y(s;x,u)\|_{\cH} \longrightarrow 0,
\ee
uniformly with respect to both the control $u$ in $\mathcal{U}_{ad}$ and the time $s$ {\mk in} $[0, T]$; see Assumption {\bf (A7)} in \cite[Section 2.4]{CKL17}.  Here,
\be
\Pi^\perp_N:=\mathrm{Id}_{\cH}-\Pi_N.
\ee
Easily checkable conditions ensuring such a double uniform convergence are identified 
in \cite[Section 2.7]{CKL17} {\mk for a broad class of evolution equations and their Galerkin approximations} but unfortunately do not apply to the case of the GK approximations considered here. The scope of this article does not allow for an investigation of such conditions in the case of DDEs, and results along this direction will be communicated elsewhere.

Nevertheless, error estimates can be still derived in the case of DDEs by adapting elements provided in \cite[Section 2.6]{CKL17}. We turn now to the formulation of such error estimates.

\subsubsection{Error estimates}

Our aim is to derive error estimates for GK approximations to the following type of optimal control problem associated with DDEs \eqref{Eq_DDE} framed into the abstract form \eqref{eq:DDE_abs0}, 
\be \label{P_sec3}  \tag {$\mathcal{P}$}
\begin{aligned}
\min \, J(x,u)  \quad \text{ s.t. }  \quad (y, u) \in L^2(0,T; \cH) \times  \mathcal{U}_{ad} \text{ solves} ~\eqref{eq:DDE_abs_c},  
\end{aligned}
\ee
in which $J$ is given by \eqref{J_tx} with $t=0$, and $y$ solves $\eqref{eq:DDE_abs_c}$ with $ y(0)  = x$ in $\cH$.

To do so, we consider the following set of assumptions

\bi
\item[(i)] The mappings $\mathcal{F}:\cH \rightarrow \cH$, $\mathcal{G}:\cH\rightarrow \mathbb{R}^+$ are locally Lipschitz,  and $\mathcal{E}:V\rightarrow \mathbb{R}^+$ is continuous.

\item[(ii)]  The {\mk linear} operator $\mathcal{A}_N: \cH_N \rightarrow \cH_N$ is diagonalizable over $\mathbb{C}$ and there exists $\alpha$ in $\mathbb{R}$ such that for each $N$
\be \label{eq_eign_bounds}
\mathrm{Re}\; \lambda_j^N \leq \alpha, \; \forall j \in \{1,\cdots,N\},
\ee
where $\{\lambda_j^N\}$ denotes the set of (complex) eigenvalues of $\mathcal{A}_N$.

\item[(iii)] Let $\mathcal{U}_{ad}$ be given by  \eqref{Eq_U_bounded} with $U$ bounded in $V$.  For each $T>0$, $(x,u)$ in $\cH \times \mathcal{U}_{ad}$, the problem \eqref{eq:DDE_abs_c} (with $t=0$) admits a unique mild solution {\mk $y(\cdot; x, u)$} in $C([0,T],\cH)$, and for each $N\geq 0$, its GK approximation \eqref{eq:DDE_Galerkin_abs} admits a unique solution $y_N(\cdot; x_N, u)$ in $C([0,T],\cH)$. Moreover, there exists a constant $\mathcal{C}:=\mathcal{C}(T,x)$ such that
\bea \label{Eq_y_uniform-in-u_bounds}
& \|y(t; x, u)\|_{\cH} \le \mathcal{C}, \qquad \Forall t\in[0,T], \; u \in  \mathcal{U}_{ad}, \\
& \|y_N(t; x_N, u)\|_{\cH} \le \mathcal{C}, \qquad \Forall t\in[0,T],\; N \in \mathbb{N}^\ast, \; u \in  \mathcal{U}_{ad}.
\eea

\item[(iv)] The mild solution $y(\cdot; x, u)$ belongs to $C([0,T],D(\mathcal{A}))$ if $x$ lives in  $D(\mathcal{A})$.
\ei

\br
\hspace*{2em}  \vspace*{-0.4em}
\bi
\item[(i)] Compared with the case of eigen-subspaces considered in \cite[Section 2.6]{CKL17}, the main difference here lies in the regularity assumption (iv) above. This assumption is used to handle the term $\Pi_N \mathcal{A} \Pi_{N}^\perp y$ arising in the error estimates; see e.g.~\eqref{energy est:1} below. Note that this term is zero when $\mathcal{A}$ is self-adjoint and $\cH_N$ is an eigen-subspace of $\mathcal{A}$, which is the setting considered in \cite[Section 2.6]{CKL17}.

\item[(ii)] A way to ensure Condition (iv) to hold consists of proving that the mild solution belongs to  $C^1([0,T],\cH)$. For conditions on $\mathcal{F}$ ensuring such a regularity see e.g.~\cite[Theorem 3.2]{Bensoussan_al07}. 
\ei
\er

Finally, we introduce the cost functional $J_N$ to be the cost functional $J$ assessed along the Galerkin approximation 
$y_N$ solving \eqref{eq:DDE_Galerkin_abs} with $ y_N(0)  = x_N \in \cH_N$, namely
\be  \label{J_Galerkin}
J_N(x_N,u) := \int_0^T [\mathcal{G}(y_N(s; x_N, u)) + \mathcal{E}(u(s)) ] \, \d s, 
\ee 
along with the optimal control problem
\be \label{P_N}  \tag {$\mathcal{P}_N$}
\begin{aligned}
\min \, J_N(x_N,u)   \;\; \text{ s.t. }  \;\; (y_N, u) \in L^2(0,T; \cH_N)  \times  \mathcal{U}_{ad} \text{ solves} ~ \eqref{eq:DDE_Galerkin_abs},
\end{aligned}
\ee
\bes
\hspace{2cm}\mbox{ with } y_N(0)  = x_N \in \cH_N.
\ees

We are now in position to derive error estimates as formulated in Theorem \ref{Thm_PM_val} below.
For that purpose Table \ref{tab_notations} provides a list of the main symbols necessary to a good reading of the proof {\mk and statement} of this theorem. The proof is largely inspired from that of \cite[Theorem 2.3]{CKL17}; the main amendments being specified within this proof.  

\begin{table}[h] 
\caption{Glossary of principal symbols}
\label{tab_notations}       
\centering
\begin{tabular}{ll}
\toprule\noalign{\smallskip}
Symbol & Terminology \\ 
\noalign{\smallskip}\hline\noalign{\smallskip}
$(y^\ast, u^\ast)$ & An optimal pair of the optimal control problem \eqref{P_sec3} \\ 
$u_{t,x}^{*}$ & Minimizer of the value function $v$ defined in \eqref{Eq_val_fcts_sec3}\\
$u^{N,\ast}_{t,x_{N}}$ & Minimizer of the value function $v_N$ defined in \eqref{Eq_val_fcts_sec3}\\
$ y_{t,x}(\cdot; u^\ast_{t,x})$ & Mild solution to \eqref{eq:DDE_abs_c} driven by $u^*_{t,x}$\\
$y_{t,x}(\cdot; u^{N,\ast}_{t,x_{N}})$ & Mild solution to \eqref{eq:DDE_abs_c} driven by $u^{N,\ast}_{t,x_{N}}$\\
$\mathfrak{B}$ & The closed ball in $\cH$ centered at $0_{\cH}$ with radius $\mathcal{C}$ given by \eqref{Eq_y_uniform-in-u_bounds}\\
$\|\mathcal{A}\|_{\rm{op}}$ & Norm of $\mathcal{A}$ as a linear bounded operator from $D(\mathcal{A})$ to $\cH$\\
\noalign{\smallskip} \bottomrule 
\end{tabular}
\end{table}

\bt\label{Thm_PM_val}
Assume that assumptions (i)-(iv) hold. Assume furthermore that for each $(t,x)$  in $[0,T)\times D(\mathcal{A})$, there exists a minimizer $u^*_{t,x}$ (resp.~$u^{N,\ast}_{t,x_{N}}$) for the value function $v$ (resp.~$v_N$) defined in \eqref{Eq_val_fcts_sec3}.

Then for any $(t,x)$  in $[0,T)\times D(\mathcal{A})$ it holds that
\bea\label{PM_value_est_goal}
 |v(t,x) - &v_N(t,x_{N})| \\
 & \le \Lip(\mathcal{G}\vert_{\mathfrak{B}}) \left[\sqrt{T-t} +  (T-t)\sqrt{f(T)} \right] \mathcal{R}(u^\ast_{t,x},u^{N,\ast}_{t,x_{N}}),
\eea
with
\be \label{eq_f}
f(T):=\bigg(\mathrm{Lip}(\mathcal{F}|_{\mathfrak{B}}) +\|\mathcal{A}\|_{\rm{op}}^2 \bigg)\exp\bigg(2\Big[\alpha + \frac{3}{2} \mathrm{Lip}(\mathcal{F}|_{\mathfrak{B}}) +\frac{1}{2}\Big]T\bigg),
\ee
and
\be
 \mathcal{R}(u^\ast_{t,x},u^{N,\ast}_{t,x_{N}}):= \|\Pi_N^\perp y_{t,x}(\cdot; u^\ast_{t,x})\|_{L^2(t,T; D(\mathcal{A}))} + \| \Pi_N^\perp y_{t,x}(\cdot; u^{N,\ast}_{t,x_{N}})\|_{L^2(t,T; D(\mathcal{A}))}.
\ee
\et

\bp
We provide here the main elements of the proof that needs to be amended from \cite[Theorem 2.3]{CKL17}. 
First note that mild solutions to \eqref{eq:DDE_abs_c} lie in $C([0,T],D(\mathcal{A}))$ due to Assumption (iv) and since the initial datum $x$ is taken in $D(\mathcal{A}).$

We want to prove that for each $N$, there exists a constant $\beta>0$ such that for any  $(t,x) \in [0,T)\times D(\mathcal{A})$, and $u\in \mathcal{U}_{ad}[t,T]$,
\be  \label{low_mode_est}
\|\Pi_{N}y_{t,x}(s;u) -y^N_{t,x_{N}}(s;u) \|^2_{\cH} \le  f(T) \int_t^s \| \Pi_N^\perp y_{t,x}(s';u) \|^2_{D(\mathcal{A})} \, \d s', \qquad  s \in [t, T].
\ee
Let us introduce 
\be\label{defw}
w(s):= \Pi_{N}y_{t,x}(s;u) -y^N_{t,x_{N}}(s;u). 
\ee 
By applying $\Pi_N$ to both sides of Eq.~\eqref{eq:DDE_abs_c}, we obtain that $\Pi_{N}y_{t,x}(\cdot; u)$  (denoted by $\Pi_{N}y$ to simplify), satisfies:
\beas
\frac{\d \Pi_{N}y}{\d s} &= \mathcal{A}_N \Pi_N y + \Pi_N \mathcal{A} \Pi_{N}^\perp y + \Pi_{N} \mathcal{F}(\Pi_N y + \Pi_{N}^\perp y) +  \Pi_{N}\mathfrak{C} (u(s)).
\eeas

This together with \eqref{eq:DDE_Galerkin_abs} implies that $w$ satisfies the following problem:
\bea \label{eq:w}
\frac{\d w}{\d s} &= \mathcal{A}_{N} w + \Pi_N \mathcal{A} \Pi_{N}^\perp y+ \Pi_{N}  \mathcal{F}(\Pi_N y + \Pi_{N}^\perp y) -  \Pi_{N}F(y_N), \\
w(t)  &= 0.
\eea
By taking the $\cH$-inner product on both sides of \eqref{eq:w} with $w$, we obtain:
\bea  \label{energy est:1}
\frac{1}{2}\frac{\d \|w\|^2_{\cH}}{\d s} = \langle \mathcal{A}_{N} w, w \rangle_{\cH} &+ \langle \Pi_N \mathcal{A} \Pi_{N}^\perp y,w\rangle_{\cH}  \\
&+\langle \Pi_{N} \bigl( \mathcal{F}(\Pi_N y + \Pi_{N}^\perp y) -  \mathcal{F}(y_N)\bigr), w \rangle_{\cH} .
\eea
We estimate now the term  $\langle \mathcal{A}_{N} w, w \rangle_{\cH}$ in the above equation using Assumption (ii). For this purpose, note that for any $p$ in $\cH_N$, the following identity holds:
\be \label{eq_AN_est1}
\langle \mathcal{A}_{N} p, p \rangle_{\cH}  = \langle M \xi, \xi \rangle_{\mathbb{R}^n},
\ee
where $M$ is the {\mk matrix} representation of $\mathcal{A}_{N}$ under the Koornwinder basis given by \eqref{eq:A}, and $\xi$  {\mk denotes the column vector} $(\xi_0, \cdots, \xi_{N-1})^{\mathrm{tr}}$ {\mk whose entries are given by}
\bes
\xi_j :=  \frac{1}{\|\mathcal{K}_j\|_{\mathcal{E}}}\langle p, \mathcal{K}^\tau_j \rangle_{\cH} , \quad j = 0, \cdots, N-1. 
\ees
Note that the $N\times N$ matrix $M$ has the same eigenvalues $\{\lambda_j^N\;:\; j = 1, N\}$ as $\mathcal{A}_N$. Thanks to Assumption (ii), $M$ is diagonalizable over $\mathbb{C}$. Denoting {\mk by $f_j$} the normalized eigenvector of $M$ corresponding to each $\lambda^N_j$, we have 
\bes
M \xi = M \sum_{j=1}^N \widetilde{\xi}_j f_j = \sum_{j=1}^N \lambda_j^N \widetilde{\xi}_j f_j,
\ees
where $\widetilde{\xi}_j = \langle \xi, f_j \rangle_{\mathbb{C}^n}$. It follows then
\bes 
\langle M \xi, \xi \rangle_{\mathbb{R}^n} = \langle \sum_{j=1}^N \lambda_j^N \widetilde{\xi}_j f_j, \xi \rangle_{\mathbb{C}^n} = \sum_{j=1}^N  \lambda_j^N |\widetilde{\xi}_j|^2 = \sum_{j=1}^N \mathrm{Re}\,\lambda_j^N |\widetilde{\xi}_j|^2,
\ees
{\mk for which the latter} equality holds since {\mk $\xi\mapsto \langle M \xi, \xi \rangle_{\mathbb{R}^n}$ is real-valued}. {\mk Due to the bound} \eqref{eq_eign_bounds}, we {\mk have thus}
\be \label{eq_AN_est2}
\langle M \xi, \xi \rangle_{\mathbb{R}^n} \le \alpha \sum_{j=1}^N |\widetilde{\xi}_j|^2.
\ee
{\mk On the other hand, by noting that}
\be \label{eq_AN_est3}
|\widetilde{\xi}_j|^2 =  \langle \widetilde{\xi}, \widetilde{\xi} \rangle_{\mathbb{C}^n} = \langle \xi, \xi \rangle_{\mathbb{R}^n} = \langle p, p \rangle. 
\ee
{\mk we obtain from \eqref{eq_AN_est1}--\eqref{eq_AN_est3}, by taking $p=w$,} that 
\be
\langle \mathcal{A}_{N} w(s), w(s) \rangle_{\cH} \leq \alpha  \|w(s)\|_{\cH}^2, \; \; s\in [t,T].
\ee

{\mk We infer from \eqref{energy est:1} that}
\beas 
\frac{1}{2}\frac{\d \|w(s)\|_{\cH}^2}{\d s} \le &\left(\alpha + \frac{3}{2} \mathrm{Lip}(\mathcal{F}|_{\mathfrak{B}}) +\frac{1}{2} \right) \|w(s)\|_{\cH}^2 \\
&+ \frac{1}{2}\|\mathcal{A}\|_{\rm{op}}^2 \| \Pi_N^\perp y\|_{D(\mathcal{A})}^2+ \frac{1}{2} \mathrm{Lip}(\mathcal{F}|_{\mathfrak{B}})  \| \Pi_N^\perp y\|_{\cH}^2,
\eeas 
in which we have also used Assumption (iii) and the local Lipschitz property of $\mathcal{F}$ on the closed ball $\mathfrak{B}$ in $\cH$ centered at $0_{\cH}$ with radius $\mathcal{C}$ given by \eqref{Eq_y_uniform-in-u_bounds}.

By Gronwall's lemma and recalling that $w(t)=0$, we obtain thus 
\bes  
 \|\Pi_{N}y_{t,x}(s;u) -y^N_{t,x_{N}}(s;u)\|_{\cH}^2  \le f(T) \int_t^s \| \Pi_N^\perp y_{t,x}(s';u)\|_{D(\mathcal{A})}^2 \d s',
\ees
with
\bes 
f(T)=\bigg(\mathrm{Lip}(\mathcal{F}|_{\mathfrak{B}}) +\|\mathcal{A}\|_{\rm{op}}^2 \bigg)\exp\bigg(2\Big[\alpha + \frac{3}{2} \mathrm{Lip}(F|_{\mathfrak{B}}) +\frac{1}{2}\Big]T\bigg). 
\ees 
Then by noting that 
\beas 
|J_{t, x}(u) -  J_{t,x_{N}}^N(u)| \le \Lip(\mathcal{G}\vert_{\mathfrak{B}}) \!\! \int_{t}^T  \!\!\Big( & \| \Pi_{N}y_{t,x}(s;u) - y^N_{t,x_N}(s;u)\|_{D(\mathcal{A})} \\
& \qquad +  \| \Pi_{N}^\perp y_{t,x}(s;u) \|_{D(\mathcal{A})}  \Big) \d s,
\eeas
we have
\bea  \label{Eq_J_est}      
& |J_{t, x}(u) -  J_{t,x_{N}}^N(u)| \\
&  \le \Lip(\mathcal{G} \vert_{\mathfrak{B}}) \left[\sqrt{T-t} +(T-t) \sqrt{f(T)} \right] \| \Pi_N^\perp y_{t,x}(\cdot;u)\|_{L^2(t,T; D(\mathcal{A}))}, 
\eea
for all $u \in  \mathcal{U}_{ad}[t,T].$

Finally by noting that
\bes 
v(t,x)  = J_{t,x}(u^*_{t,x}) \le J_{t,x}(u^{N,\ast}_{t,x_{N}}) \quad \text{and} \quad  v_N(t,x_{N}) = J^N_{t,x_{N}}(u^{N,\ast}_{t,x_{N}}), 
\ees
we obtain
\bes 
 v(t,x) - v_N(t,x_{N})  \le J_{t,x}(u^{N,\ast}_{t,x_{N}}) - J^N_{t,x_{N}}(u^{N,\ast}_{t,x_{N}}). 
 \ees
It follows then from \eqref{Eq_J_est} that
\bea \label{PM_val_est1}
&  v(t,x)  - v_N(t,x_{N}) \\
 & \;\; \le  \Lip(\mathcal{G} \vert_{\mathfrak{B}}) \left[{\sqrt{T-t}}+ (T-t) \sqrt{f(T)}\right] 
 \| \Pi_N^\perp y_{t,x}(\cdot;u^{N,\ast}_{t,x_{N}})\|_{L^2(t,T; D(\mathcal{A}))}.
\eea
Similarly, 
\bea \label{PM_val_est2}
& v_N(t,x_{N})  -  v(t,x) \\
& \;\; \le  \Lip(\mathcal{G}\vert_{\mathfrak{B}}) \left[{\sqrt{T-t}}+ (T-t)\sqrt{f(T)} \right] 
 \|\Pi_N^\perp y_{t,x}(\cdot; u^\ast_{t,x})\|_{L^2(t,T; D(\mathcal{A}))}.
\eea
The estimate \eqref{PM_value_est_goal} results then from \eqref{PM_val_est1} and \eqref{PM_val_est2}.

\ep

{\mk We conclude this section with the following corollary providing the error estimates between the optimal control and that obtained from a  GK approximation.}

\bc\label{Lem_controller_est}
Assume that the conditions given in Theorem~\ref{Thm_PM_val} hold. {\mk Given $x$ in $\cH$}, let us {\mk introduce the notations}, $u^* := u^*_{0,x}$ and $u^{*}_{N}:= u^{N,*}_{0,x_N}$. Assume furthermore that there exists $\sigma >0$ such that the following local growth condition is satisfied  for the cost functional $J$ (with $t=0$) given by \eqref{J_tx}: 
\be\label{Eq_growth_onJ}
\sigma \|u^*  - v\|_{L^q(0,T; V)}^q \le J(x, v) - J(x, u^*), 
\ee
for all $v$ in some neighborhood $\mathcal{O} \subset \mathcal{U}_{ad}$ of $u^*$, with $\mathcal{U}_{ad}$ given by \eqref{Eq_U_bounded}. Assume finally that $u^*_N$  lies in  $\mathcal{O}$. 
Then, 
\bea\label{Est_contr_diff}
& \|u^\ast - u^\ast_{N}\|_{L^q(0,T; V)}^q  \\
& \le \frac{1}{\sigma}\Lip(\mathcal{G}\vert_{\mathfrak{B}}) \left[\sqrt{T} +T\sqrt{f(T)} \right] \Bigl( \| \Pi_N^\perp y^\ast(\cdot;u^*)\|_{L^2(0,T; D(\mathcal{A}))} \\
&\hspace{6cm} + 2 \|  \Pi_N^\perp y (\cdot; u^{*}_{N})\|_{L^2(0,T; D(\mathcal{A}))} \Bigr), 
\eea
where the constant $f(T)$ is given by \eqref{eq_f}.
\ec

\bp
By the assumptions, we have
\bes
\|u^*  - u^*_N\|_{L^q(0,T; V)}^q \le \frac{1}{\sigma} \left(J(x, u^*_N) - J(x, u^*)\right).
\ees
Note also that
\beas
J(x, u^*_N) - J(x, u^*) &= J(x, u^*_N) - J_N(x_{N},u^*_N) + J_N(x_{N}, u^*_N) - J(x, u^*) \\
& =  J_{0,x}(u^*_N) - J^N_{0,x_{N}}(u^*_N)+ v_N(0,x_{N}) - v(0,x),
\eeas
where we used the fact that 
\beas
 & J(x, u^*_N) = J_{0,x}(u^*_N),  \quad J_N(x_{N},u^*_N) = J^N_{0,x_{N}}(u^*_N),  \\
 &J_N(x_{N},u^*_N)  = v_N(0,x_{N}) \quad \text{and} \quad   J(x, u^*) = v(0,x).
\eeas
The result follows by applying the estimate \eqref{Eq_J_est} to $J_{0,x}(u^*_N) - J^N_{0,x_{N}}(u^*_N)$ and the estimate \eqref{PM_value_est_goal} to $v_N(0,x_{N}) - v(0,x)$.

\ep

\section{Application to an oscillation amplitude reduction problem} \label{Sect_DDE}

\subsection{Optimal control of the Wright equation} \label{Sect_Wright_PMP}
{\mk As an application,} we consider the Wright equation
\be \label{eq:Wright}
\frac{\d m}{\d t} =  - m(t-\tau) (1 +  m(t)),
\ee
where $m$ is the unknown  scalar function, and $\tau$ is a nonnegative delay parameter. 
{\mk This equation has been studied by numerous authors, and notably among the first ones, are Jones \cite{jones1962existence,jones1962existence2}, Kakutani and Markus \cite{Kakutani1958},  and Wright \cite{Wright55}.}
This equation can also be transformed via a simple change of variable \cite{Ruan06} into the well-known Hutchinson equation \cite{Hutchinson48} arising in population dynamics. {\mk It corresponds then to} the logistic equation {\mk with a delay effect introduced into the}  {\it intraspecific competition term}, namely with the change of variable $m=-1+p$, Eq.~\eqref{eq:Wright} becomes
\be\label{Eq_logistic}
\frac{\d p}{\d t} =  p(t) (1 -p(t-\tau) ).
\ee
Essentially the idea of Hutchinson \cite{Hutchinson48} was to point out that negative effects that high  population density $p$ has on the environment, influences birth rates at later times due to developmental and maturation delays, justifying thus Eq.~\eqref{Eq_logistic}.

As a result, the solutions of Eq.~\eqref{eq:Wright} are obtained from a simple shift of the solutions of Eq.~\eqref{Eq_logistic}, and reciprocally. Since Eq.~\eqref{Eq_logistic} benefits of a more intuitive interpretation, we will often prefer to think about Eq.~\eqref{eq:Wright} in terms of Eq.~\eqref{Eq_logistic}.  

{\mk Anyway, it is known that Eq.~\eqref{eq:Wright} (and thus Eq.~\eqref{Eq_logistic})} experiences a supercritical Hopf bifurcation {\mk as the delay parameter is varied and crosses a critical value, $\tau_c$,} from below; {\mk see e.g.~\cite[Sect.~9]{Chow_al77} and \cite{sun2007analysis}}.\footnote{The critical value $\tau_c$ at which the trivial steady state of the Wright equation \eqref{eq:Wright} loses its linear stability is given by $\tau_c = \frac{\pi}{2}$. This can be found by analyzing the associated linear eigenvalue problem $-\phi(\theta - \tau) = \beta \phi(\theta), \theta \in [-\tau, 0]$. By using the ansatz $\phi(\theta) = e^{\beta \theta}$, we obtain $-e^{-\beta \tau} = \beta$. Assuming  that $\beta = \pm i\omega $, we get $-\cos(\omega \tau) + i \sin(\omega \tau) = i \omega $, leading thus to $\omega = 1$ and $\tau = \frac{(2n+ 1)\pi}{2}$, $n = 0, 1, 2, \cdots$. Consequently,  the critical delay parameter is $\tau_c = \frac{\pi}{2}$.} 

\begin{figure}[ht]
\begin{center}
    \includegraphics[keepaspectratio=true,scale=0.34]{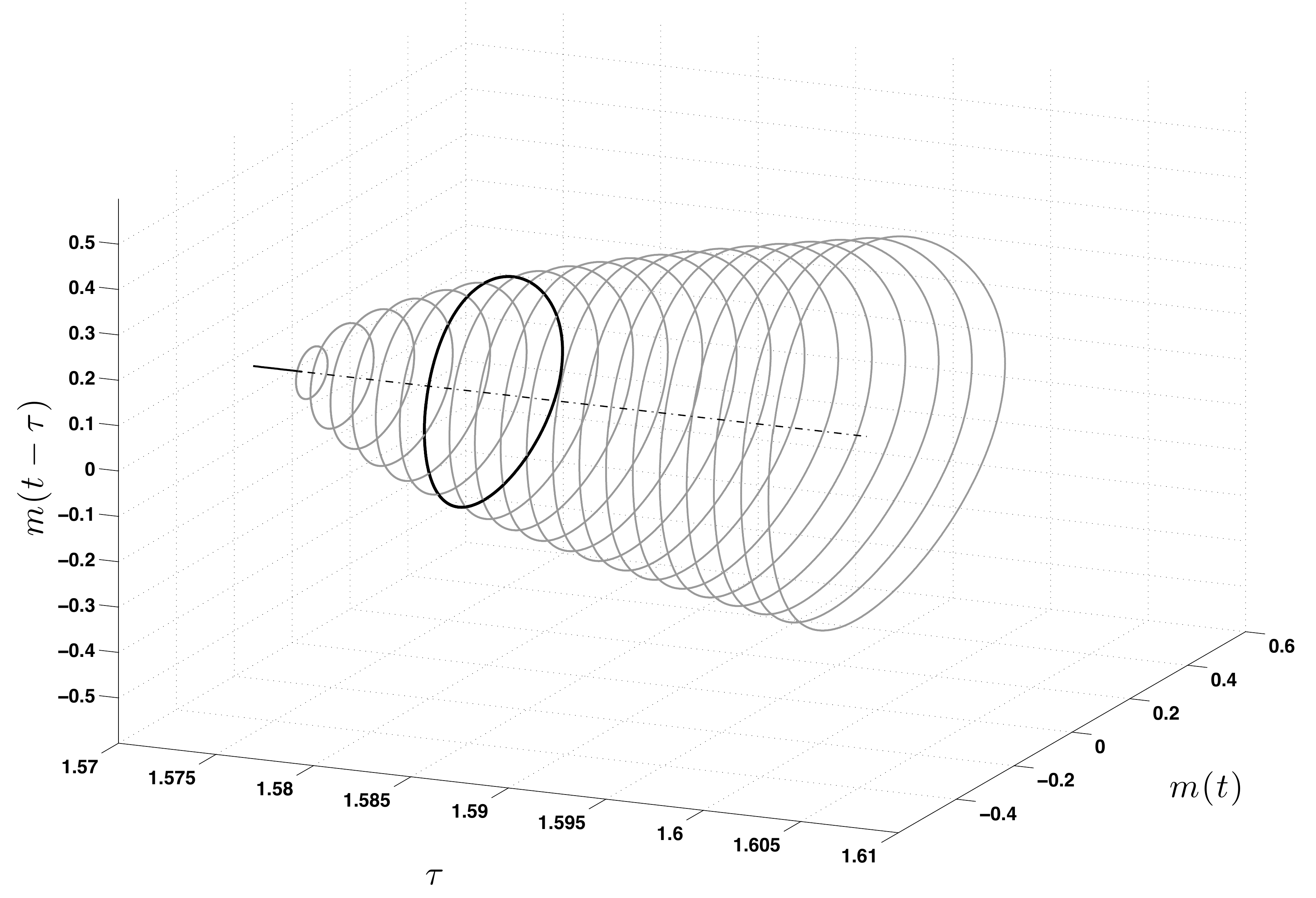}
   \caption{{\footnotesize {\mk {\bf Unfolding of the Hopf bifurcation taking place for Eq.~\eqref{eq:Wright} as $\tau$ is varied}. The limit cycles are represented in the  the embedded reduced phase space $(m(t),m(t-\tau))$. The limit cycle in heavy black indicates the  ``$\tau$-secion'' considered for the numerical results of Sectns.~\ref{Sec_num_results} and \ref{Sec_HJB} below. The trivial steady state is represented via a dotted line when unstable and a plain one when stable.}}}
\label{Fig_Hopf}
\end{center}
\end{figure}

{\mk Figure \ref{Fig_Hopf} shows---in the embedded reduced phase space $(m(t),m(t-\tau))$---the corresponding bifurcation of the limit cycle unfolding from the trivial steady state, as $\tau$ is varied\footnote{Such an embedded phase space is classically used to visualize attractors associated with DDEs; see e.g.~\cite{CGN17}.}. As Fig.~\ref{Fig_Hopf} illustrates, the amplitude of the oscillation that takes place via the Hopf bifurcation, increases as $\tau$ is increased away from $\tau_c$.}

{\mk Since fluctuating populations are susceptible to extinction due to sudden and unforeseen environmental disturbances \cite{RC11}, a knowledge of the conditions for which the population is resilient to these disturbances  is of great interest in planning and designing control as well as management strategies; see e.g.~\cite{CR06,RC07,roques2010does}.}

Given a convex functional depending on $m$ and a control to be designed, our goal is to show {\mk (for Eq.~\eqref{eq:Wright})} that {\mk relatively} close to the criticality (i.e.~$\tau \approx \tau_c$), the {\mk oscillation amplitude} can be reduced at a nearly optimal cost, by solving efficiently low-dimensional optimal control problems.

{\mk For that purpose, we consider, given  $u$ in $L^2(0,T; \mathbb{R})$ and $T > 0$, the following forced version of Eq.~\eqref{eq:Wright}} 
\be \label{eq:Wright_forced}
 \frac{\d m}{\d t} =  - m(t-\tau) (1 +  m(t)) + u(t), \qquad t \in (0,T), 
\ee
supplemented by the initial history  
$ m(t) = \phi(t)$  (with $\phi \in L^2(-\tau,0;\mathbb{R})$) for $t$ in the interval $[-\tau,0)$, and  by $m(0) = m_0$ (in $\mathbb{R}$) at $t=0$. 
{\mk In Eq.~\eqref{eq:Wright_forced}, $u(t)$ can be thought as an environmental management strategy. Our goal is to understand the management strategies that may lead to a reduction of the oscillation  amplitude of the population (and possibly a stabilization); i.e.~to determine $u$ for which  $m\approx 0$ and thus $p\approx1$, while dealing with limited resources to plan such strategies.}

Naturally, we introduce thus the following cost functional
\be  \label{J_DDE}
J(m,v) := \int_0^T \left[\frac{1}{2}m(t)^2 + \frac{\mu}{2} u(t)^2 \right] \d t, \; \mu > 0.
\ee 
We are interested in solving optimal control problem associated with Eq.~\eqref{eq:Wright_forced} for this cost functional.  
To address this problem, we recast {\mk first} this optimal control problem into the abstract framework of Sect.~\ref{Sec_DDEreformulation}.

{\mk In that respect}, we introduce  the variables 
\be\label{state-variable}
y:= (m_t, m_t(0)), \; \mbox{ with } m \mbox{ solving } \eqref{eq:Wright_forced},
\ee
and
we take $\mathcal{U}_{ad}$ defined in \eqref{Eq_U_bounded} with $q=2$ and  with $U$ to be (possibly) a bounded subset of $V:=\mathbb{R}$.

The control operator $\mathfrak{C}\colon V \rightarrow \mathcal{H}$ appearing in Eq.~\eqref{eq:DDE_abs_c}  is taken here to be linear and given by
\bea \label{Def_C2}
[\mathfrak{C} v ](\theta) & := \begin{cases}
0, &  \theta \in[-\tau, 0)   \vspace{0.4em}\\ 
v, & \theta = 0
\end{cases},\; \;  v \in V.
\eea

As explained in Sect.~\ref{Sec_DDEreformulation},
Eq.~\eqref{eq:Wright_forced} can be recast into the following abstract ODE posed on $\cH$
\be \label{eq:DDE_abs_ex}
\frac{\d y}{\d t}   = \mathcal{A} y + \mathcal{F}(y) +  \mathfrak{C} u(t), 
\ee
with  the linear operator $\mathcal{A}$ given in \eqref{Def_A} that, in the case of  Eq.~\eqref{eq:Wright_forced}, takes the form 
\bea \label{Def_A_ex}
\lbrack \mathcal{A} \Psi \rbrack (\theta) & := \begin{cases}
{\displaystyle \frac{\d^+ \Psi^D}{\d \theta}}, &  \theta \in[-\tau, 0),  \vspace{0.4em}\\ 
{\displaystyle -\Psi^D(-\tau)}, & \theta = 0,
\end{cases} 
\eea
with domain given in \eqref{D_of_A}. 
The nonlinearity $\mathcal{F}$ given in \eqref{Def_F} takes here the form, 
\bea \label{Def_F_ex}
[\mathcal{F} (\Psi) ](\theta) & := \begin{cases}
0, &  \theta \in[-\tau, 0),   \vspace{0.4em}\\ 
-\Psi^D(-\tau)\Psi^S, & \theta = 0, 
\end{cases}  \quad \Forall \Psi = (\Psi^D, \Psi^S) \in  \mathcal{H}.
\eea
Let $\Phi$ be in $\cH$ {\mk given by}
\be \label{eq_Wright_initial_abs}
\Phi:= (\phi, m_0), \; \mbox{ with } \phi \in L^2(-\tau,0;\mathbb{R}), \; m_0 \in \mathbb{R}. 
\ee

We rewrite now the cost functional defined in \eqref{J_DDE} by using the state-variable $y$ given in \eqref{state-variable}, and 
define thus
\be  \label{J_DDE_recast}
J(y,u) := \int_0^T \left[\frac{1}{2}\big(m_t(0)\big)^2 + \frac{\mu}{2} u(t)^2 \right] \d t;
\ee 
still denoted by $J$.

The optimal control problem associated with Eq.~\eqref{eq:DDE_abs_ex} becomes then
\be \label{P1_recast} 
\min \, J(\Phi, u) \; \text{s.t.} \; (y,u) \in L^2(0,T; \cH) \times \mathcal{U}_{ad} \text{ solves } \text{Eq.~\eqref{eq:DDE_abs_ex} with } y(0)=\Phi \in \cH,
\ee
with $\mathcal{U}_{ad}$ defined in \eqref{Eq_U_bounded} with $q=2$ and with either $U=\mathbb{R}$ or $U$ to be a bounded subset of $V=\mathbb{R}$.

\br
Note that in the case $u(t)$ in Eq.~\eqref{eq:Wright_forced} is replaced by $u(t-r)$ (with $0<r$), one can still
recast Eq.~\eqref{eq:Wright_forced} into an abstract ODE of the form \eqref{eq:DDE_abs_ex}; i.e.~to deal with control with delay.
In that case  by taking $V:=H^1(-r,0;\mathbb{R})\times \mathbb{R}$, the operator $\mathfrak{C}$  becomes  
\bea \label{Def_C3}
[\mathfrak{C} (\phi^D,\phi^S)](\theta) & := \begin{cases}
0, &  \theta \in[-\tau, 0)   \vspace{0.4em}\\ 
\phi^D(-r), & \theta = 0
\end{cases},\; \;  (\phi^D,\phi^S) \in V.
\eea
This reformulation is compatible with the framework of \cite{CGLW15} and thus allow us to derive GK approximations of such abstract ODEs. We leave the associated optimal control problems to the interested reader. 
\er

\subsection{Galerkin-Koornwinder approximation of the optimal control problem} \label{Sect_Wright_GK_control}

{\mk First, we specify for Eq.~\eqref{eq:Wright_forced}, the (concrete) ODE version \eqref{Galerkin_AnalForm} of the (abstract) $N^{\mathrm{th}}$ GK approximation \eqref{eq:DDE_Galerkin_abs}.} For this purpose, we first compute $\Pi_N \mathfrak{C} u(s)$.  

Note that since $\{\mathcal{K}^\tau_j \; : \; j \in \mathbb{N}\}$  forms a Hilbert basis of $\cH$, for any $\Psi=(\Psi^D,\Psi^S)$ in $\mathcal{H}$, we have (see also \cite[Eq.~(3.25)]{CGLW15})
\bea\label{Eq_decomp}
\Psi &=\sum_{j=0}^{\infty} \frac{\langle \Psi,\mathcal{K}^\tau_j \rangle_{\mathcal{H}}}{\|\mathcal{K}^\tau_j\|_{\mathcal{H}}^2}\mathcal{K}_j^\tau\\
&  =\sum_{j=0}^{\infty} \Big(\frac{1}{\tau}\langle \Psi^D,K_j^\tau \rangle_{L^2}+\Psi^S K_j^\tau(0)\Big) \frac{\mathcal{K}_j^\tau }{\|\mathcal{K}_j^\tau \|_{\mathcal{H}}^2}.
\eea

By taking now $\Psi=\mathfrak{C} v$ with $\mathfrak{C}$ given by \eqref{Def_C2} and $v$ in $\mathbb{R}$, we have
\be
\mathfrak{C} v=\sum_{j=0}^{\infty} \Big(v K_j^\tau(0)\Big) \frac{\mathcal{K}_j^\tau }{\|\mathcal{K}_j^\tau \|_{\mathcal{H}}^2},
\ee
which leads to 
\be
\Pi_N \mathfrak{C} v=v \sum_{j=0}^{N-1}   \frac{\mathcal{K}_j^\tau }{\|\mathcal{K}_j^\tau \|_{\mathcal{H}}^2},
\ee
recalling that $K_j^\tau(0)=1$; see \eqref{Eq_normalization}.

Then by defining  $\mathfrak{C}_N \colon V \rightarrow \mathbb{R}^N$ to be
\be \label{eq:C_N}
\mathfrak{C}_N v := \Bigl ( \frac{1}{\|\mathcal{K}_0 \|_{\mathcal{E}}^2}, \cdots, \frac{1}{\|\mathcal{K}_{N-1} \|_{\mathcal{E}}^2 }  \Bigr )^{\tr} v, 
\ee
one can write the Galerkin approximation \eqref{eq:DDE_Galerkin_abs} in the Koornwinder  basis as the ODE system
\bea \label{eq:Wright_Galerkin}
& \frac{\d \boldsymbol{\xi}_N}{\d t}  = M \boldsymbol{\xi}_N + G (\boldsymbol{\xi}_N) + \mathfrak{C}_N u(t),\\
& \boldsymbol{\xi}_N(0)  = \boldsymbol{\zeta}_N \in \mathbb{R}^N.
\eea
Here $\boldsymbol{\xi}_N = (\xi^N_0, \cdots, \xi^N_{N-1})^{\mathrm{tr}}$; the $N\times N$-matrix $M$ given in \eqref{eq:A} has its elements given here by
\begin{equation} \label{eq:A2}
\begin{aligned}
(M)_{j,n} & = \frac{1}{\|\mathcal{K}_j\|_{\mathcal{E}}^2 }\Big(- K_n(-1) + \frac{2}{\tau}\sum_{k=0}^{n-1} a_{n,k} \left( \delta_{j,k} \|\mathcal{K}_j\|^2_{\mathcal{E}} - 1 \right ) \Big), \; j, n = 0, \cdots, N-1;
\end{aligned}
\end{equation}
Recall that here the Kronecker symbol $\delta_{j,k}$ has been used, and  the coefficients $a_{n,k}$ are obtained by solving a triangular linear system in which the right hand side has 
explicit coefficients depending on $n$; see \cite[Prop.~5.1]{CGLW15}.

{\mk The initial datum} $\boldsymbol{\zeta}_N:= (\zeta^N_0, \cdots, \zeta^N_{N-1})^{\mathrm{tr}}$ has its components defined by
\be \label{eq:Wright_Galerkin_IC} 
\zeta^N_j = \frac{1}{\|\mathcal{K}_{j}\|_{\mathcal{E}}^2}  \langle \Phi, \mathcal{K}_j^\tau \rangle_{\mathcal{H}}, \;\; \; j = 0, \cdots, N-1, 
\ee
where $\Phi$ in $\mathcal{H}$ denotes the initial data {\mk given by \eqref{eq_Wright_initial_abs}} for the abstract ODE \eqref{eq:DDE_abs_ex}.

The nonlinearlity $G$ given in \eqref{eq:G} takes, in the case of {\mk Eq.~\eqref{eq:Wright}}, the form  
\be \label{eq:G2}
G_j(\boldsymbol{\xi}_N) =  - \frac{1}{\|\mathcal{K}_{j}\|_{\mathcal{E}}^2} \Bigg[\sum_{n=0}^{N-1} \xi^N_n(t) \Bigg] \Bigg[\sum_{n=0}^{N-1} \xi^N_n(t)K_n(-1) \Bigg], \; 0\leq j\leq N-1.
\ee

{\mk 
We turn next to the rewritting---using the variables $\xi^N_n$---of the corresponding cost functional $J_N$  given in its abstract version by \eqref{J_Galerkin}. To do so, it is sufficient to recall from \eqref{eq_mN} that the $N$-dim GK approximation $m_N$ is given by}
\bes
m_N(t) = \sum_{n=0}^{N-1} \xi^N_{n}(t) K_n^\tau(0) = \sum_{n=0}^{N-1} \xi^N_n(t),
\ees
{\mk which, denoting still by $J_N$ the rewritten cost functional, gives} 
\bea \label{cost_JN}
J_N(\boldsymbol{\zeta}_N,u) := \int_0^T \left[\frac{1}{2}\Big(\sum_{n=0}^{N-1} \xi^N_n(t; \boldsymbol{\zeta_N},u) \Big)^2 + \frac{\mu}{2} u(t)^2 \right] \d t, \; \boldsymbol{\zeta}_N \in \mathbb{R}^N.
\eea

{\mk We are now in position to write the corresponding Galerkin approximations of the initial optimal control problem \eqref{P1_recast}. 
More precisely, the (reduced) optimal control problem associated with the $N$-dimensional GK approximation \eqref{eq:Wright_Galerkin} of Eq.~\eqref{eq:Wright_forced}, is}
\be  \label{PG} 
\hspace*{-0.5ex}\boxed{
\begin{aligned}
 & \min \, J_N(\boldsymbol{\zeta}_N, u)  \text{ subject to } (\boldsymbol{\xi}_N(\cdot; \boldsymbol{\zeta}_N, u), u) \in L^2(0,T; \mathbb{R}^N) \times  {\mk \mathcal{U}_{ad}}  \\
&\text{ solves the {\mk ODE system} \eqref{eq:Wright_Galerkin} supplemented with $\boldsymbol{\xi}_N(0) = \boldsymbol{\zeta}_N$,}
\end{aligned}
}
\ee
{\mk with $\mathcal{U}_{ad}$ defined in \eqref{Eq_U_bounded} with $q=2$ and with either $U=\mathbb{R}$ or $U$ to be a bounded subset of $V=\mathbb{R}$; see Sect.~\ref{Sec_HJB}.}

\subsection{Numerical results}\label{Sec_num_results}

As mentioned earlier, the delay parameter $\tau$ in Eq.~\eqref{eq:Wright} is taken to be slightly above its critical value $\tau_c$ at which the supercritical Hopf bifurcation takes place; see again Figure~\ref{Fig_Hopf}. In particular, the $\tau$-value is selected so that there is only one conjugate pair of unstable eigenmodes for the linearized DDE. We can thus hope for relatively low-dimensional Galerkin systems aimed at the synthesis of (sub)optimal controls at a nearly optimal cost \cite[Sect.~4]{CKL17}, and obtained in a feedback form from the associated reduced HJB equation if the reduced dimension is low enough; see 
Sect.~\ref{Sec_HJB} below.

Indeed, for this choice of $\tau$-value, a very good approximation of the optimal control can be already synthesized from the 6D-GK system via a PMP approach, as checked by comparison with higher-dimensional GK systems.  To further reduce the dimension, we will show that a 2D projected GK system obtained from the 6D-GK system allows for the synthesis of a sub-optimal control at a nearly optimal cost (see Table~\ref{table_errors}) whereas, the 2D-GK system fails to do so. 

We describe next the precise setup of the optimal control problem. After that, we provide the explicit forms of the 2D- and 6D-GK systems, and explain how the 2D projected GK system is obtained. The control synthesis results by using the PMP approach are then presented at the end of this section. As a comparison, the results by using the HJB approach are presented in Sect.~\ref{Sec_HJB}.

\subsubsection{The numerical setup}

For the delay parameter, we take $\tau=1.58$. Although this $\tau$ value represents less than $1\%$ increase from the critical value $\tau_c = \pi/2$, the amplitude of the emerged stable periodic oscillation is already not so small, as can be observed from Figure~\ref{Fig_Hopf} (the thick black limit cycle). The parameter $\mu$ in the cost functional $J$ given by \eqref{J_DDE_recast} is taken to be $1/2$. For the selection of initial history segment and the time horizon of the control, the uncontrolled Wright equation \eqref{eq:Wright} is integrated over a sufficiently long time interval $[0, t^*]$ so that the solution has already settled down to the attracting periodic orbit. The control horizon $T$ is chosen such that a half of the oscillation period has developed\footnote{In our experiments it corresponds to $T=4$.} for the uncontrolled problem (see the dashed curve in Figure~\ref{Fig:DDE_PMP_result}), and we take then a ``snippet'' of this periodic orbit over $[t^\ast-T-\tau, t^\ast-T]$ as the initial history for the controlled problem \eqref{eq:Wright_forced}.  In other words, the DDE \eqref{eq:Wright_forced} is initiated from a snippet of the attracting periodic orbit, $\Phi$, from which we want to reduce its amplitude at an optimal cost by solving \eqref{P1_recast}. We turn hereafter to the determination of the relevant GK systems in order to propose an approximate solution to this problem.

\subsubsection{The concrete Galerkin-Koornwinder approximation systems}

Recall that the $N$-dim GK system of Eq.~\eqref{eq:Wright_forced} is given by Eq.~\eqref{eq:Wright_Galerkin}--\eqref{eq:G2}. To get access to the numerical values of the coefficients involved in Eq.~\eqref{eq:Wright_Galerkin}, we only need to specify the matrix $M$, the square of the norm of the basis function $\|\mathcal{K}_j \|_{\mathcal{E}}^2$ as well as the value $K_j(-1)$ of the Koornwinder polynomials at the left endpoint for $j=0,\cdots, N-1$. 

In our case, we have (cf.~\cite[Proposition 3.1]{CGLW15})
\bea
\big(\|\mathcal{K}_{0}\|_{\mathcal{E}}^2, \; \cdots, \|\mathcal{K}_{5}\|_{\mathcal{E}}^2 \big) &  =  (2,  \;  3.333,   \; 10, \; 24.2857,   \; 49.1111,   \; 87.4545), \\
(K_0(-1), \; \cdots, K_5(-1)) & = (1, \;  -3,  \;   7,  \;  -13,  \;  21, \;  -31).
\eea
With $N=6$, we have from \eqref{eq:A2} that the matrix $M$ appearing in Eq.~\eqref{eq:Wright_Galerkin} is given after computation by 
\be \label{eq_M6}
M_{6} = \begin{pmatrix}
-0.5000  &  2.7658  & -5.3987  & 10.9304  & -16.8291  &  25.6266 \\
-0.3000  &  0.1405  &  2.8367  &  -2.5557  &   6.6114   & -7.4089 \\
-0.1000  &  0.0468  & -2.2190  &   8.6418  & -9.8215   &  18.4165 \\
-0.0412  &  0.0193  & -0.9137  &  -1.6538  &  8.4652   &  -7.0109 \\
-0.0204  &  0.0095  & -0.4518  &  -0.8178  & -3.2628   & 11.8690 \\
-0.0114  &  0.0054  & -0.2537  &  -0.4593  & -1.8323   & -3.1193 
\end{pmatrix},
\ee
up to four digits precision after the decimal place. 

The initial data for the 6D-GK system is simply taken to be the projection of the aforementioned DDE initial history segment $\Phi$ onto the first 6 Koornwinder basis functions; cf.~\eqref{eq:Wright_Galerkin_IC}. We have thus
\be \label{eq_6D_initial}
\boldsymbol{\zeta}_6 = (0.0590,\;   0.0827,\;   0.0014,\;   -0.0006,\;   0.0,\;    0.0)^{\tr}.
\ee 

For the 2D-GK sytem, the matrix $M_2$ consists of the $2\times2$ block in the  upper left corner of $M_6$ given by \eqref{eq_M6}. Namely,
\be \label{eq_M2}
M_{2} = \begin{pmatrix}
-0.5000  &  2.7658   \\
-0.3000  &  0.1405 
\end{pmatrix}.
\ee
The initial data for the 2D-GK system is the same as the first two component of $\boldsymbol{\zeta}_6$ given by \eqref{eq_6D_initial}. 

As reported in Table~\ref{table_errors} below, the 6D-GK system produces a control strategy that is already optimal as compared with the cost from the 12D-GK system. On the other hand, the 2D-GK system leads to a control strategy that significantly inflates the control cost. This failure of the 2D-GK system can be understood by simply inspecting its linear part. 

Indeed, although the first two Koornwinder modes already capture more than $98\%$ of the $L^2$-energy contained in the uncontrolled DDE solution, the linear part of the corresponding 2D-GK system does not resolve well the pair of unstable eigenvalues of the DDE's linear part. More precisely, while the eigenvalues of $M_2$ are given by 
\be
\beta_{2,1} = -0.1798 + 0.8527 \,i, \quad  \beta_{2,2} = -0.1798 - 0.8527 \,i,
\ee
the first pair of eigenvalues of the DDE's linear part is given up to four significant digits by
\be\label{Eq_osc_part}
\beta_{1} = 0.0026 + 0.9958 \, i, \quad  \beta_{2} = 0.0026 - 0.9958 \, i.
\ee 
This unstable pair of eigenvalues can nevertheless be resolved up to the given precision from $M_6$ given in \eqref{eq_M6}, which actually corresponds to
the lowest dimension to do so.  

Although, the 2D-GK system lacks to resolve important dynamical features such as the dominant eigenpair of the DDE's linear part, we provide for later usage the explicit expression of the 2D-GK system: 
\be \label{eq_2D_GK}
 \frac{\d \boldsymbol{\xi}_2}{\d t} = M_2 \boldsymbol{\xi}_2  - \begin{pmatrix}
 0.5 (\xi_{2,1})^2  - \xi_{2,1}\xi_{2,2} - 1.5 (\xi_{2,2})^2 \vspace*{0.2ex}\\
0.3 (\xi_{2,1})^2 - 0.6 \xi_{2,1}\xi_{2,2} - 0.9 (\xi_{2,2})^2
 \end{pmatrix} + 
 \begin{pmatrix}
  0.5 \\
  0.3
  \end{pmatrix} u(t),
\ee
where $\boldsymbol{\xi}_2 := (\xi_{2,1}, \; \xi_{2,2})^{\tr}$, and the initial data is $\boldsymbol{\xi}_2(0) =(0.0590,\;   0.0827)^{\tr}$. 

Let $J_2$ be the cost functional given by \eqref{cost_JN} with $N=2$, the control problem of the form \eqref{PG} associated with \eqref{eq_2D_GK} reads
\be  \label{P_2D} 
\hspace*{-0.5ex}\boxed{
\begin{aligned}
 & \min \, J_2(\boldsymbol{\zeta}_2, u)  \text{ subject to } (\boldsymbol{\xi}_2(\cdot; \boldsymbol{\zeta}_2, u), u) \in L^2(0,T; \mathbb{R}^2) \times \mathcal{U}_{ad} \\
&\text{ solves the ODE system \eqref{eq_2D_GK} supplemented with $\boldsymbol{\xi}_2(0) = \boldsymbol{\zeta}_2$,}
\end{aligned}
}
\ee
where the initial data $\boldsymbol{\zeta}_2$ is taken to be $(0.0590,\;   0.0827)^{\tr}$.

\subsubsection{The projected GK system}

In this section we propose a natural way to cope with the resolution issue of the dominant eigenpair \eqref{Eq_osc_part}, while keeping the dimension 
of the reduced system as small as possible. It is based on the projection onto the space $\widehat{\mathcal{H}}_2$ spanned by the first two eigenvectors\footnote{Note that the eigenvalues are labeled in the lexicographical order.  Namely,  for $1\le n < n' \le 6$, we have either $\mathrm{Re} \beta_n > \mathrm{Re} \beta_{n'}$, or  $\mathrm{Re} \beta_n = \mathrm{Re} \beta_{n'}$ and $\mathrm{Im} \beta_n \ge  \mathrm{Im} \beta_{n'}$.} of the matrix $M_6$ given by \eqref{eq_M6}; recalling that the sixth dimension constituting for this example the minimal dimension to resolve accurately \eqref{Eq_osc_part}.

Note that since the eigenvalues of $M_6$ are all complex-valued, the resulting 2D projected GK system has complex-valued coefficients. However, the eigenvectors of $M_6$ appear in conjugate pairs because $M_6$ is a real-valued matrix. As a result, the two components of the projected system, denoted by $\boldsymbol{z}_2:=(z_{2,1},\, z_{2,2})^{\tr}$, are complex conjugate to each other. We can thus rewrite the system under a new real-valued variable $\widetilde{\boldsymbol{\xi}}_2:=(\widetilde{\xi}_{2,1}, \, \widetilde{\xi}_{2,2})^{\tr}$ defined by
\be \label{eq_cplx_transform}
\widetilde{\xi}_{2,1} =  \frac{z_{2,1} + z_{2,2}}{2}, \qquad \widetilde{\xi}_{2,2} =  \frac{z_{2,1} - z_{2,2}}{2 \, i}. 
\ee

Using the {\sc Matlab} built-in function \texttt{eig} to compute the eigenbasis of $M_6$ and the {\sc Matlab} Symbolic Math Toolbox, the resulting 
transformed equations for the $\widetilde{\boldsymbol{\xi}}_2$-variable are given by
\be \label{eq_2D-proj}
 \frac{\d \widetilde{\boldsymbol{\xi}}_2}{\d t} \!= \! \begin{pmatrix}
 0.0026 & -0.9958  \\
  0.9958   & 0.0026
\end{pmatrix} \widetilde{\boldsymbol{\xi}}_2  + \begin{pmatrix}
  b^1_{20} (\widetilde{\xi}_{2,1})^2 +   b^1_{11} \widetilde{\xi}_{2,1} \widetilde{\xi}_{2,2} + b^1_{02}(\widetilde{\xi}_{2,2})^2 \vspace*{0.2ex}\\
b^2_{20} (\widetilde{\xi}_{2,1})^2 +   b^2_{11} \widetilde{\xi}_{2,1}\widetilde{\xi}_{2,2} + b^2_{02}(\widetilde{\xi}_{2,2})^2
 \end{pmatrix} + 
 \begin{pmatrix}
  \alpha_1 \\
  \alpha_2
  \end{pmatrix} u(t),
\ee
with the coefficients given by $\alpha_1 = 0.0608$, $\alpha_2 = 0.1133$, and
\beas 
& b_{2,0}^1 = -1.7887, && b_{1,1}^1 =  2.1915, &&  b_{0,2}^1 = 1.7996, \\
& b_{2,0}^2 = -3.3320, && b_{1,1}^2 = 4.0824, &&  b_{0,2}^2 = 3.3524.
\eeas
The initial data is given by 
\be \label{eq_2D-proj-initial}
\widetilde{\boldsymbol{\xi}}_2(0) = (0.0107 , 0.0253)^{\mathrm{tr}},
\ee
obtained by projecting the initial data of the 6D Galerkin system given by \eqref{eq_6D_initial} onto $\widehat{\mathcal{H}}_2$,  followed by the transformation \eqref{eq_cplx_transform}. 

Associated with the cost functional $J_6$ defined by \eqref{cost_JN} for the $6$D GK system, we have the following cost functional $\widetilde{J}_2$ associated with the 2D projected GK system \eqref{eq_2D-proj}:
\be
\widetilde{J}_2(\widetilde{\boldsymbol{\zeta}}_2, u):=\int_0^T \left[\frac{1}{2} \Big( d_1 \widetilde{\xi}_{2,1}(t; \widetilde{\boldsymbol{\zeta}}_2, u) + d_2 \widetilde{\xi}_{2,1}(t; \widetilde{\boldsymbol{\zeta}}_2, u)  \Big)^2 + \frac{\mu}{2} u(t)^2 \right]\d t, \; \widetilde{\boldsymbol{\zeta}}_2 \in \mathbb{R}^2,
\ee
with $d_1 = -4.0668$ and $d_2 = 7.2710$. 

Then, corresponding to the optimal control problem \eqref{PG}, the optimal control problem associated with \eqref{eq_2D-proj}--\eqref{eq_2D-proj-initial} reads: 
\be  \label{P_2D-proj} 
\hspace*{-0.5ex}\boxed{
\begin{aligned}
 & \min \, \widetilde{J}_2(\widetilde{\boldsymbol{\zeta}}_2, u)  \text{ subject to } (\widetilde{\boldsymbol{\xi}}_2(\cdot; \widetilde{\boldsymbol{\zeta}}_2, u), u) \in L^2(0,T; \mathbb{R}^2) \times  \mathcal{U}_{ad}  \\
&\text{ solves the ODE system \eqref{eq_2D-proj} supplemented with $\widetilde{\boldsymbol{\xi}}_2(0) = \widetilde{\boldsymbol{\zeta}}_2$,}
\end{aligned}
}
\ee
where the initial data $\widetilde{\boldsymbol{\zeta}}_2$ is given by \eqref{eq_2D-proj-initial}.

\subsubsection{Numerical results obtained from the PMP}

We present now the numerical results obtained by applying the PMP to the the optimal control problems associated with the three ODE systems specified above, namely, the problem \eqref{PG} associated with the GK system for $N=2$ and $N=6$, and the problem \eqref{P_2D-proj} associated with the 2D projected GK system \eqref{eq_2D-proj}.

\begin{figure}[ht]
   \centering
   \includegraphics[height=0.4\textwidth, width=1\textwidth]{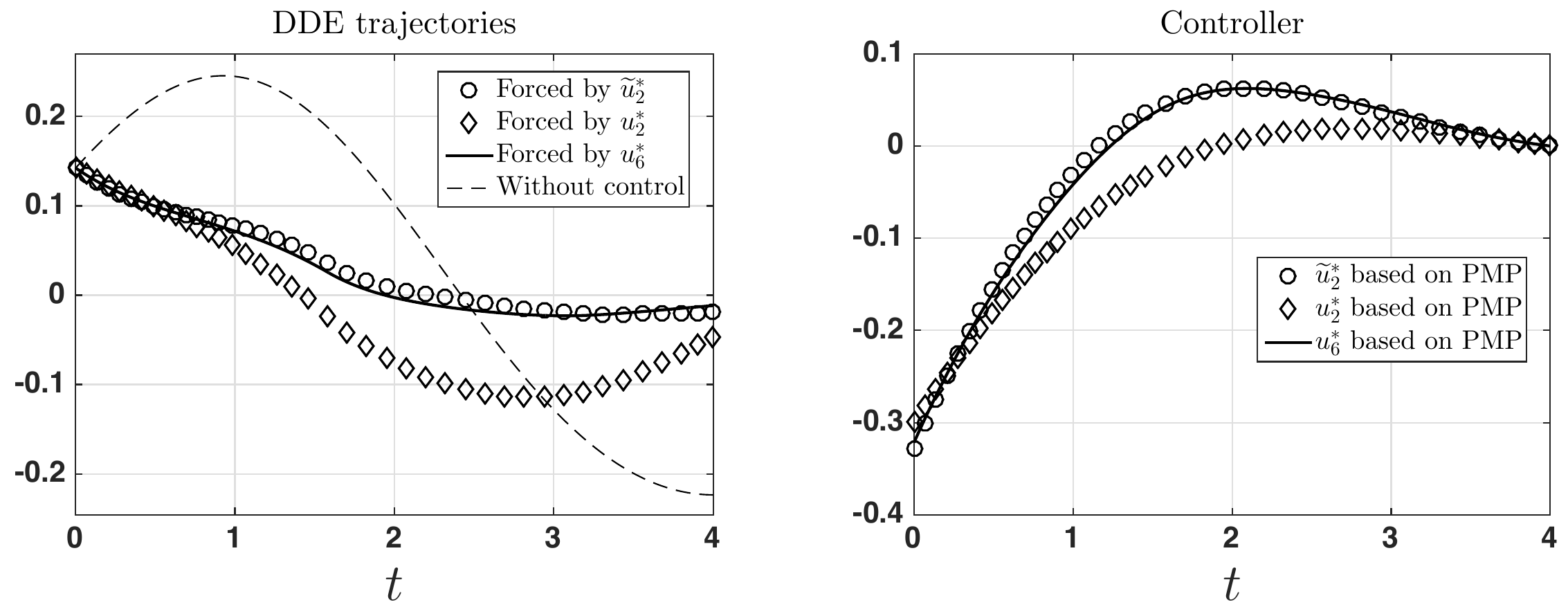}
\vspace{-2ex}
  \caption{{\bf Left panel:} DDE trajectories by solving Eq.~\eqref{eq:Wright_forced} driven respectively by the (sub)optimal control $\widetilde{u}^*_2$, $u^*_2$, and $u^*_6$; the uncontrolled trajectory is also plotted as a comparison. {\bf Right panel:} The controller $\widetilde{u}^*_2$ synthesized by \eqref{P_2D-proj} based on the 2D projected GK system, and the controllers $u^*_2$ and $u^*_6$ synthesized by \eqref{PG} based respectively on the 2D- and 6D-GK systems. We have taken $\mu = 0.5$ in the cost functional \eqref{J_DDE}. The maximal time of control is $T = 4$ and the delay parameter is taken to be $\tau= 1.58$.}   \label{Fig:DDE_PMP_result}
\end{figure}

In each case, the PMP approach leads to a boundary value problem (BVP) for the corresponding state and co-state variables, which is solved using the {\sc Matlab} built-in function \texttt{bvp4c}; see e.g.~\cite[Sect.~5]{CL15} for more details. The results are shown in Fig.~\ref{Fig:DDE_PMP_result}. As already mentioned earlier, the 6D-GK system allows for a control synthesis that is nearly optimal by comparison with that synthesized from the 12D-GK system; cf.~Table~\ref{table_errors}. In contrast, the control synthesized from the 2D-GK system leads to a substantially higher cost value.
Remarkably, the 2D projected GK system \eqref{eq_2D-proj}  allows for a control synthesis at a cost value very close to that obtained from higher-dimensional GK systems, i.e.~at a nearly optimal cost. This success is due simply to the fact that the latter system resolves the unstable pair of eigenvalues of the DDE's linear part, while the former does not.

\begin{table}[ht]
\caption{Cost values and relative errors}
\label{table_errors}       
\centering
\begin{tabular}{ccccc}
\toprule\noalign{\smallskip}
 & 12D-GK  & 6D-GK &  2D-GK & 2D projected GK \\
 \noalign{\smallskip}\hline\noalign{\smallskip}
 Cost values & 0.0163 & 0.0163 & 0.0253 & 0.0166 \\
 Relative error & & &  54.93\% & 1.4655\%  \\
\noalign{\smallskip} \bottomrule 
\end{tabular}
\end{table}

\section{Towards the approximation of HJB equations associated with DDEs}\label{Sec_HJB}

 In this section we compute optimal controls by solving the HJB equations associated with the control problems  \eqref{P_2D} and \eqref{P_2D-proj}. This allows to obtain (numerically) the controls in feedback form. However, the approach becomes infeasible in higher dimensions by the curse of dimensions (i.e. the exponential increase of complexity with the number of unknowns), since the dimension of the underlying Galerkin system is equal to the dimension of the space the HJB equation is defined on (with respect to the spatial variable).

Let us introduce for $(\boldsymbol{\eta},u) \in \R^2 \times V$ the functionals
\begin{subequations} \label{l_funct}
\begin{align}
\mathcal{L}(\boldsymbol{\eta},u)&:=\frac{1}{2}|\boldsymbol{\eta}|^2 + \frac{\mu}{2} u^2,\label{case1} \\  
 \widetilde{\mathcal{L}}(\boldsymbol{\eta},u)&:=\frac{1}{2}( -4.0668 \eta_1 + 7.2710 \eta_2)^2 + \frac{\mu}{2} u^2.\label{case2}
 \end{align}
 \end{subequations}
 Here the functional $\mathcal{L}$ is associated with  the 2D-GK system  \eqref{eq_2D_GK}, while the functional  $ \widetilde{\mathcal{L}}$ is associated with the 2D projected GK system \eqref{eq_2D-proj}.

Given a control $u$ in $\mathcal{U}_{ad}[t,T]$, we define then the following cost functionals over the interval $[t,T]$ 
\begin{align}
J_{t,\boldsymbol{\eta}}(u):=\int_t^T \mathcal{L}(\boldsymbol{\xi}_2(s;\boldsymbol{\eta},u),u(s))\d s, \quad \widetilde{J}_{t,\boldsymbol{\eta}}(u):=\int_t^T  \widetilde{\mathcal{L}}(\widetilde{\boldsymbol{\xi}}_2(s;\boldsymbol{\eta},u),u(s))\d s
\end{align}
where $\boldsymbol{\xi}_2(s;\boldsymbol{\eta},u)$ denotes the solution of the  2D-GK system \eqref{eq_2D_GK},  and $\widetilde{\boldsymbol{\xi}}_2(s;\boldsymbol{\eta},u)$ the solution of the 2D projected GK system  \eqref{eq_2D-proj}; solutions emanating in each case from the initial datum $\boldsymbol{\eta}$ in $\R^2$.
 The associated value functions are given by
\begin{align}
v(t,\boldsymbol{\eta}):=\inf_{u\in \mathcal U_{ad}[t,T]}J_{t,\boldsymbol{\eta}}(u), \quad \widetilde{v}(t,\boldsymbol{\eta}):=\inf_{u\in \mathcal U_{ad}[t,T]} \widetilde{J}_{t,\boldsymbol{\eta}}(u).
\end{align}

 Formally, the value functions are associated with the instationary HJB equation of the following type
\begin{align} \label{hjb_ex}
 \partial_t v(t,{\boldsymbol{\eta}})+H({\boldsymbol{\eta}},\nabla_{\boldsymbol{\eta}} v(t,\boldsymbol{\eta}))=0\quad \text{in } (0,T) \times \R^2,\quad v(T,{\boldsymbol{\eta}})=0 \quad \text{in } \R^2;
\end{align}
with Hamiltonian 
\begin{align}
H({\boldsymbol{\eta}},{\boldsymbol{p}}):=\inf_{u\in U} \left(\langle f({\boldsymbol{\eta}},u), {\boldsymbol{p}}\rangle_{\mathbb{R}^2} + \mathcal{L}({\boldsymbol{\eta}},u)\right),\quad ({\boldsymbol{\eta}},{\boldsymbol{p}}) \in \R^2 \times \R^2;
\end{align}
where $f$ denotes either the RHS of \eqref{eq_2D_GK} if $\mathcal{L}$ is given by \eqref{case1}, or the RHS of \eqref{eq_2D-proj} if $\mathcal{L}=\widetilde{\mathcal{L}}$ given by \eqref{case2}. To be more specific, we denote by ${\boldsymbol{\alpha}}$ the vector of coefficients in front of $u(t)$ in either Eq.~\eqref{eq_2D_GK} or Eq.~\eqref{eq_2D-proj}.  We can then write $f({\boldsymbol{\eta}},u)$ under the following generic compact form
\be
f({\boldsymbol{\eta}},u):=K {\boldsymbol{\eta}} + N({\boldsymbol{\eta}}) +{\boldsymbol{\alpha}} u,
\ee
with  $K=M_2$ and $N({\boldsymbol{\eta}})$ corresponding to the quadratic terms in Eq.~\eqref{eq_2D_GK}, or with  $K{\boldsymbol{\eta}}$ and $N({\boldsymbol{\eta}})$ corresponding respectively to the linear and  quadratic terms in Eq.~\eqref{eq_2D-proj}.

In each case, by noting that the terms depending on $u$ in $\langle f({\boldsymbol{\eta}},u), {\boldsymbol{p}} \rangle_{\mathbb{R}^2} + \mathcal{L}({\boldsymbol{\eta}},u)$,
can be rewritten as 
\bea\label{Explicit_H0}
 \frac{\mu}{2} u^2 + \langle \boldsymbol{\alpha}, \boldsymbol{p} \rangle_{\mathbb{R}^2} u & =  \frac{\mu}{2} \Big(u^2 + \frac{2}{\mu} \langle \boldsymbol{\alpha}, \boldsymbol{p} \rangle_{\mathbb{R}^2} u \Big) \\
& =  \frac{\mu}{2} \Big( \big (u + \frac{1}{\mu} \langle \boldsymbol{\alpha}, \boldsymbol{p} \rangle_{\mathbb{R}^2} \big)^2 - \frac{1}{\mu^2} \langle \boldsymbol{\alpha}, \boldsymbol{p} \rangle_{\mathbb{R}^2}^2 \Big). 
\eea
we obtain, for $U=\mathbb{R}$, that 
\be\label{Explicit_H}
H({\boldsymbol{\eta}},{\boldsymbol{p}})=\frac{1}{2}\langle {\boldsymbol{\eta}}, Q {\boldsymbol{\eta}} \rangle_{\mathbb{R}^2}  + \langle K {\boldsymbol{\eta}} + N({\boldsymbol{\eta}}),{\boldsymbol{p}} \rangle_{\mathbb{R}^2} -\frac{\langle{\boldsymbol{\alpha}},{\boldsymbol{p}} \rangle_{\mathbb{R}^2}^2}{2 \mu},
\ee
with $Q=\mbox{Id}_2$ in the case of Eq.~\eqref{eq_2D_GK}, and 
$\langle {\boldsymbol{\eta}}, Q {\boldsymbol{\eta}} \rangle_{\mathbb{R}^2}  =( -4.0668 \eta_1 + 7.2710 \eta_2)^2$ in the case of Eq.~\eqref{eq_2D-proj}.

To prove that the value functions are the unique viscosity solutions of
\eqref{hjb_ex} goes beyond the scope of this paper; here, we formulate the HJB equation
and the feedback law formally; we also assume that the value function is smooth. For an analysis of HJB equations associated with infinite
horizon optimal control problems for differential equations with distributed delays we refer to \cite{Federico_al10}.

\subsubsection{Discretization}
For solving the HJB equation \eqref{hjb_ex} there exist various schemes.
In the presented setting we are in particular interested in efficient solvers for high dimensional HJB equations to allow also the treatment of reduced systems of more than two dimensions. 
Recently, different approaches have been considered to address this challenge arising from the curse of dimensions when solving HJB equations; see e.g.~\cite{KaliseKunisch:2017} for a polynomial approximation of high-dimensional HJB equations. In the context of suboptimal control of PDEs an approach to solve the associated HJB equation based on sparse grids was considered in \cite{GarckeKroener:2016}. Here, we use a finite difference scheme method; more precisely, an {\it essentially non-oscillatory} (ENO) scheme \cite{OsherShu:1991} for space discretization is coupled with a Runge-Kutta time discretization scheme of second order following \cite{BokanowskiForcadelZidani:2010,OsherShu:1991}.

We give a brief description of the scheme for a given continuous Hamiltonian $H \colon \R^2\times \R^2\rightarrow \R$.
For temporal mesh parameter $\Delta t>0$ and spatial mesh parameter $h=(h_1,h_2)\in \R^2_{>0}$, we define a spatial mesh
 \begin{align}
 \mathcal T:=\{(k h_1, l h_2)\in \R^2 :  (k,l)\in \Z^2\}
 \end{align}
and temporal mesh 
\be
t=t_0<t_1< ...  < t_N= T,\quad \Delta t=t_{j+1}-t_j,
\ee
for  $j=0,...,N-1$ and $N$ in $\N^\ast$.
Given an approximation of the value function on the mesh $\mathcal T$ at time step $t_j$ denoted by $\boldsymbol{v}\colon \Z^2 \rightarrow \R$, where $\boldsymbol{v}_{kl}$ for $(k,l) \in \Z^2$ gives an approximation of the value function at grid point  $\boldsymbol{\eta}_{kl}:=(kh_1,lh_2 )$, we define the difference quotients 
\be
(D^+\boldsymbol{v})_{kl}:=((D^+_1 \boldsymbol{v})_{kl},(D^+_2 \boldsymbol{v})_{kl}),\quad (D^-\boldsymbol{v})_{kl}:=((D^-_1 \boldsymbol{v})_{kl},(D^-_2 \boldsymbol{v})_{kl} )
\ee
with
\bea
(D^{\pm}_1 \boldsymbol{v})_{kl}:=\pm \frac{(\boldsymbol{v}_{k\pm1,l}- \boldsymbol{v}_{kl})}{h_1},\quad
 (D^{\pm}_2 \boldsymbol{v})_{kl}:=\pm \frac{(\boldsymbol{v}_{k,l\pm 1}- \boldsymbol{v}_{kl})}{h_2}.
\eea
 Then the Lax-Friedrichs scheme
for the HJB equation reads as follows at mesh node~$\boldsymbol{\eta}_{kl}$ by: 
\begin{equation}
 \left\{
\begin{aligned}
 {\boldsymbol{v}}_{kl}^{N}&=0, && \\
 {\boldsymbol{v}}_{kl}^{j-1}&={\boldsymbol{v}}_{kl}^j-\Delta t H_{LF}({\boldsymbol{\eta}}_{kl},(D^+ {\boldsymbol{v}}^j)_{kl},(D^- {\boldsymbol{v}}^j)_{kl})
 && \text{for } j=N,\cdots,1,
\end{aligned}
\right.
\end{equation}
with $H$ given in \eqref{Explicit_H} and the Lax-Friedrichs Hamiltonian $ \Hlf \colon \R^2\times \R^2 \rightarrow \R$ given by 
\be\label{LF_Hamiltonian}
\Hlf ({\boldsymbol{\hat \eta}},\boldsymbol{p}^+,\boldsymbol{p}^-):=H\left( {\boldsymbol{{\hat \eta}}},\frac{\boldsymbol{p}^+ +\boldsymbol{p}^-}{2}\right) -\sum_{i=1}^2\frac{\nu_i({\boldsymbol{{\hat \eta}}})}{2}(p_i^+ - p_i^-),
\ee
and the stabilizing functions $\nu_i({\boldsymbol{{\hat \eta}}})$ satisfying 
\be
\max_{\boldsymbol{p} \in \R^{2}} \left|\partial_{\boldsymbol{p}_i}  H({\boldsymbol{{\hat \eta}}},\boldsymbol{p})\right| \le \nu_i({\boldsymbol{{\hat \eta}}}),
\ee
for all ${\boldsymbol{{\hat \eta}}}\in \R^2$ and $i=1$ to $2$. The Lax-Friedrichs scheme is monotone and the convergence is of first order as far as the following Courant-Friedrichs-Lewy (CFL) condition is satisfied,
\be\label{Eq_CFL}
\Delta t\max_{{\boldsymbol{{\hat \eta}}}\in \mathcal T}\left( \frac{\nu_1({\boldsymbol{{\hat \eta}}})}{h_1}+ \frac{\nu_2({\boldsymbol{{\hat \eta}}})}{h_2}\right)\le 1.
\ee
 An ENO scheme can be obtained by considering a variant of the Lax-Friedrichs scheme, namely 
\be\label{ENO_eq}
\boldsymbol{v}_{kl}^{j-1}=\boldsymbol{v}_{kl}^j-\Delta t H_{LF}({\boldsymbol{\eta}}_{kl},(\widetilde{D}^+ \boldsymbol{v}^j)_{kl},(\widetilde{D}^- \boldsymbol{v}^j)_{kl}),
\ee
where $\widetilde{D}^{\pm} \boldsymbol{v}^j$ are higher order approximations of the gradient of the value function at time step $t_j$ coupled with a
Runge-Kutta time discretization scheme of second order,
see \cite{BokanowskiForcadelZidani:2010} and \cite{OsherShu:1991}.  

The HJB equation is posed on the full space $\R^2$ for the spatial variable. For the numerical computation we restrict it to a bounded domain 
\be\label{comp_domain}
D:=[c_1,d_1] \times [c_2,d_2],
\ee
for some $\boldsymbol{c},\boldsymbol{d}$  in $\R^2$,
 such that $D$ contains the trajectory obtained from the PMP approach in Sect.~\ref{Sec_num_results}.
Choosing numbers of grid points $p_1$ and $p_2$ for each of the spatial dimension,  we set
\be
h_i:=\frac{d_i - c_i}{p_i-1},\quad i=1,2.
\ee
and obtain the spatial meshing 
\begin{align}
 \{(c_1+kh_1, c_{2}+lh_2): (k,l)\in [0,\cdots,p_1-1] \times [0,\cdots,p_2-1]\}.
\end{align}
 An appropriate choice for the boundary condition is required to get a well defined numerical scheme. 
At time steps $t_j$ for $j$ from $N$ to $1$ and grid point $\boldsymbol{\eta}_{kl}$ we define the upwind derivatives on the boundary for all $l \in [0,\cdots,p_2-1]$ by
\be
    \begin{aligned}\label{Bd1}
  k=0:&&  (D^+_{\eta_1}\boldsymbol{v})_{kl}&= \frac{\boldsymbol{v}_{k+1,l} -\boldsymbol{v}_{kl}}{h_1},\\
   k=p_1-1: && (D^-_{\eta_1}\boldsymbol{v})_{kl}&= -\frac{\boldsymbol{v}_{k-1.l} -\boldsymbol{v}_{kl}}{h_1},
    \end{aligned}
    \ee
    and for all $k \in [0,\cdots,p_1-1]$ by
    \be
    \begin{aligned}\label{Bd2}
  l=0:&&  (D^+_{\eta_2}\boldsymbol{v})_{kl}&= \frac{\boldsymbol{v}_{k,l+1}-\boldsymbol{v}_{kl}}{h_2},\\
   l=p_2-1: && (D^-_{\eta_2}\boldsymbol{v})_{kl}&= -\frac{\boldsymbol{v}_{k,l-1} -\boldsymbol{v}_{kl}}{h_2}.
    \end{aligned}
    \ee
This choice of boundary condition  corresponds to a vanishing second order derivative boundary condition. 
   
\subsubsection{Computation of optimal trajectories}\label{sec:Trajectory}
For the computation of optimal trajectories we formulate a feedback law based on the solution of the HJB equation \eqref{hjb_ex}; the solution is assumed to be smooth here. 
    Following e.g.~\cite[p.~11]{BardiCapuzzoDolcetta:2008} 
    a control $u^*$ is optimal for the initial condition $\boldsymbol{\eta}$ if and only if 
    \be
    u^*(s)=S(s,{\boldsymbol{\xi}}^*(s;\boldsymbol{\eta},u)) \quad\text{for a.e. } s\in(0,T),
    \ee
     with 
    \begin{align}\label{feedback-S}
S(s,\boldsymbol{\hat \eta }) := \operatorname{argmin}_{u\in U} 
\left(\langle f(\boldsymbol{\hat \eta},u) \nabla_{\boldsymbol{\eta}} v(s,\boldsymbol{\hat \eta})\rangle_{\mathbb{R}^2}  +\mathcal{L}(\boldsymbol{\hat \eta},u)\right)\quad  \text{for all }\boldsymbol{\hat \eta} \in \R^2,
    \end{align}
 where, here again,    $f$ denotes either the RHS of \eqref{eq_2D_GK} if $\mathcal{L}$ is given by \eqref{case1}, or the RHS of \eqref{eq_2D-proj} if $\mathcal{L}=\widetilde{\mathcal{L}}$ given by \eqref{case2}. 
To obtain an optimal trajectory one then solves the system
\begin{equation}\label{eq:ode-with-feedback}
\begin{aligned}
 \frac{\d {\boldsymbol{\xi}}(s;\boldsymbol{\eta},u)}{\d s} &= f({\boldsymbol{\xi}}(s;\boldsymbol{\eta},u),S(s,{\boldsymbol{\xi}}(s;\boldsymbol{\eta},u))),\quad s>0,\quad 
 {\boldsymbol{\xi}}(0;\boldsymbol{\eta},u)=\boldsymbol{\eta}.
\end{aligned}
\end{equation}

From \eqref{Explicit_H0}, the element $u^\ast$ that leads after minimization to the Hamiltonian $H$ given in \eqref{Explicit_H} (and thus to \eqref{feedback-S}), is given by 
\be
u^\ast = - \frac{1}{\mu} \langle \boldsymbol{\alpha}, \boldsymbol{p} \rangle_{\mathbb{R}^2}.
\ee
which leads thus to the following {\it feedback law}:
\begin{align}\label{Eq_feedback-law}
 S(s,\boldsymbol{\hat \eta})= - \frac{1}{\mu}  \langle \boldsymbol{\alpha}, \nabla_{\boldsymbol{ \eta}} v(s,\boldsymbol{\hat  \eta}) \rangle_{\mathbb{R}^2}.
 \end{align}
For the numerical computations, the gradient $\nabla_{\boldsymbol{\eta}} v$ is  approximated  component-by-component by a central difference quotient. 
      Once the feedback operator is defined, the optimal trajectory
 is derived from \eqref{eq:ode-with-feedback}.

\subsubsection{Numerical results}

Here, we compute optimal controls for \eqref{P_2D} and \eqref{P_2D-proj} by solving the corresponding HJB equation as described above. We choose the domain $D$ given in \eqref{comp_domain} in such a way that it contains the optimal trajectory by using the numerical results obtained from the PMP approach. 

For \eqref{P_2D} we set $c_1=c_2=-0.02$ and $d_1=d_2=0.1$, the temporal mesh parameter to $\Delta t = 1.543 \cdot 10^{-4}$ and the spatial one to $h_1=h_2=8.5 \cdot 10^{-3}$. The CFL constants in \eqref{Eq_CFL} are chosen to be $\nu_1=5$ and $\nu_2=2$.

For \eqref{P_2D-proj}, we set $c_1=c_2=-0.04$ and $d_1=d_2=0.04$, $\Delta t = 1.543 \cdot 10^{-4}$, and $h_1=h_2=1.33 \cdot 10^{-3}$. The CFL constants are here $\nu_1=1.5 $ and $\nu_2= 1.5$.

For the computation of the optimal control from \eqref{Eq_feedback-law} for  \eqref{P_2D}  or \eqref{P_2D-proj}, we use the corresponding spatial and temporal mesh parameters given above.

\begin{figure}[ht]
   \centering
 \includegraphics[height=0.45\linewidth, width=0.95\linewidth]{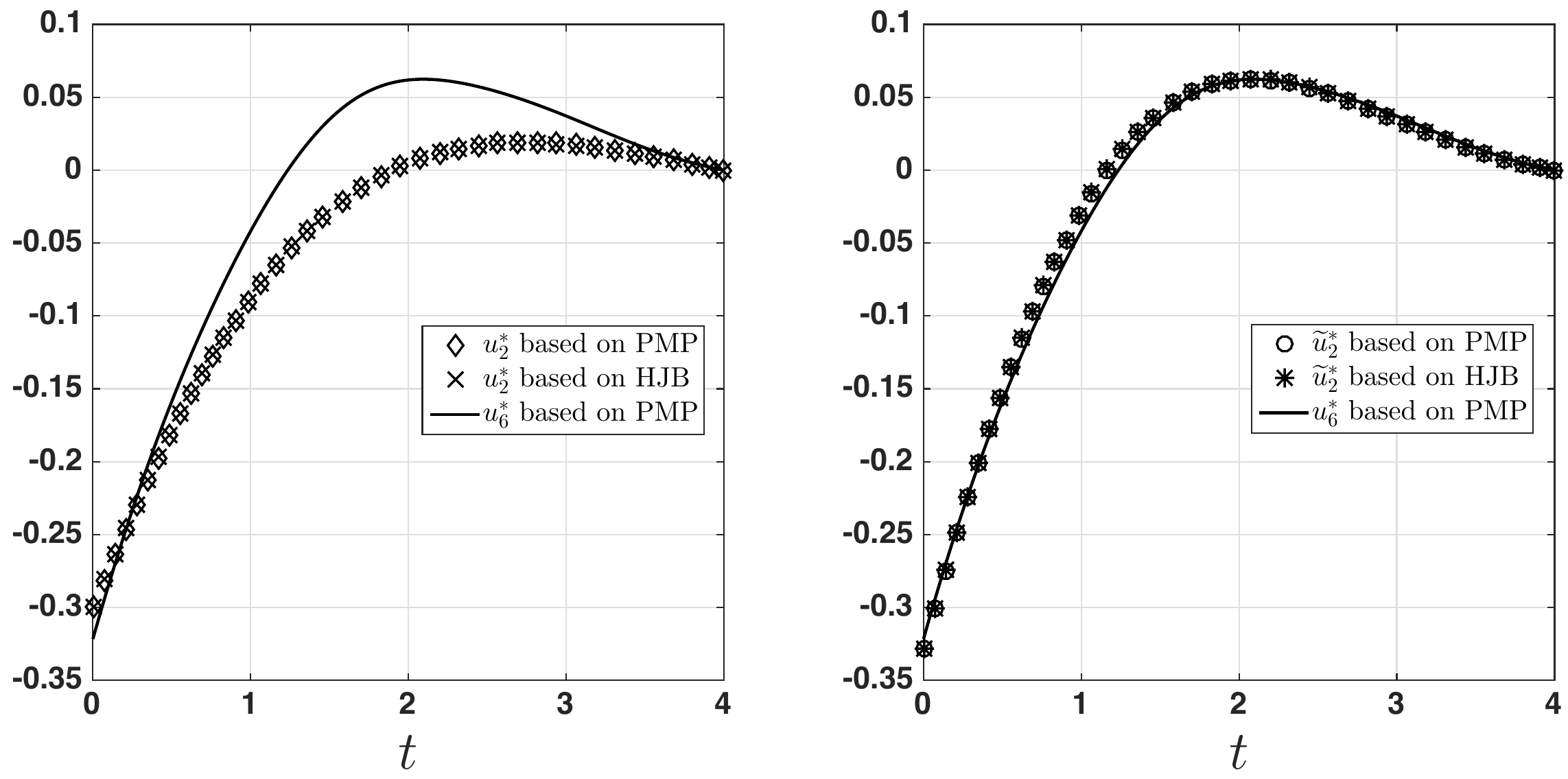}
  \caption{{\bf Left panel}: Control obtained without basis transformation (HJB vs PMP), i.e.~by solving \eqref{P_2D}. {\bf Right panel}: Control obtained after basis transformation (HJB vs PMP), i.e.~by solving \eqref{P_2D-proj}.}   \label{fig:h1}
\end{figure}

Figure \ref{fig:h1} 
shows the optimal controls for problems \eqref{P_2D} and \eqref{P_2D-proj} derived from the corresponding HJB equations as described above, and also computed from the PMP approach; see Sect.~\ref{Sec_num_results}. In each case, the syntheses of (sub)optimal controls obtained either from application of the PMP or by solving the associated HJB equations, provide nearly identical numerical solutions.

\begin{figure}[ht]
\centering
   \includegraphics[height=0.5\linewidth,, width=0.85\linewidth]{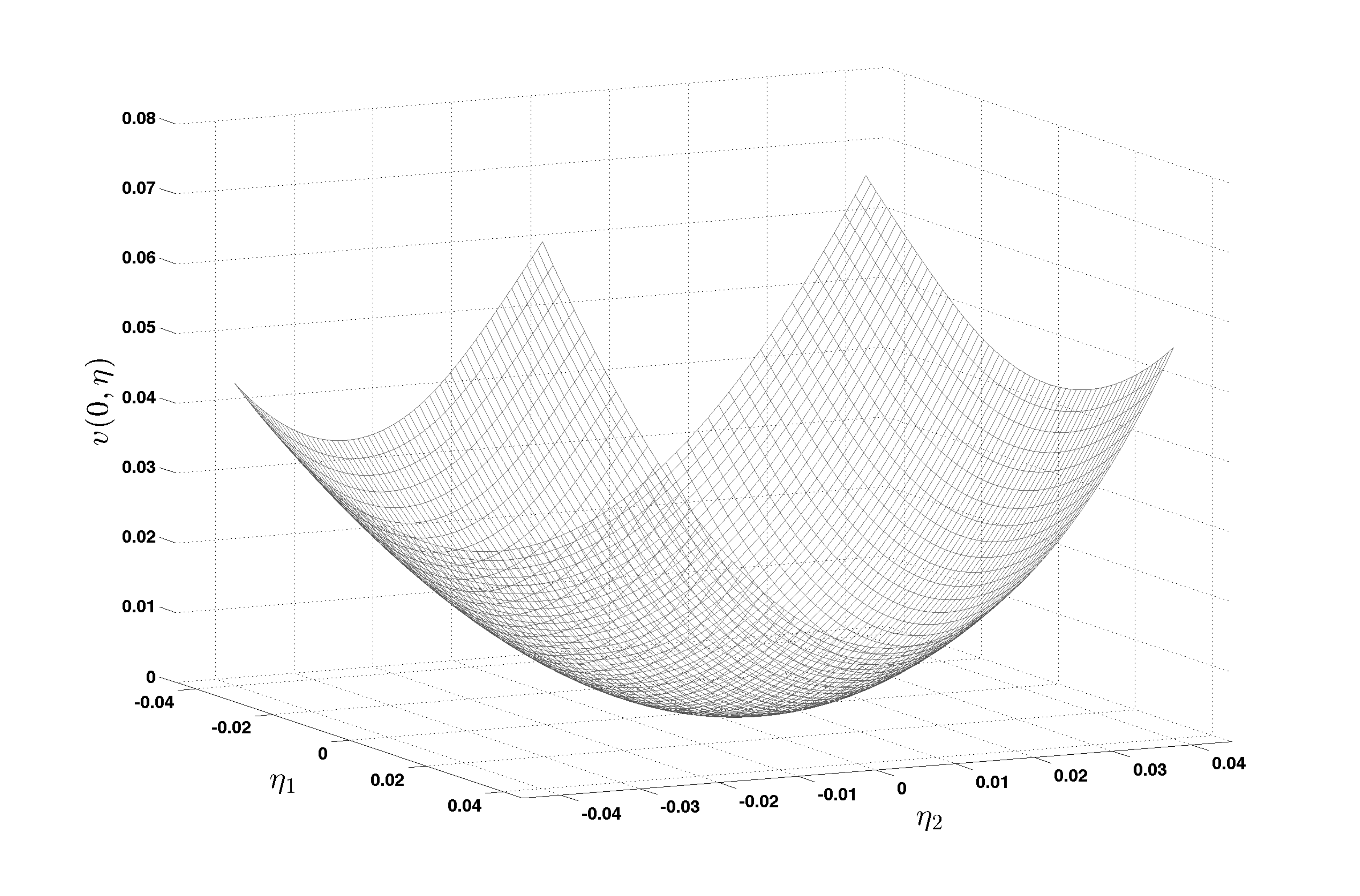}
  \caption{Value function solving, at $t=0$, the reduced HJB equation \eqref{hjb_ex} associated with \eqref{P_2D-proj}.}\label{fig:VF}
  \vspace{-2ex}
\end{figure}

Figure \ref{fig:VF} shows an approximation\footnote{According to the ENO scheme \eqref{ENO_eq} over the domain \eqref{comp_domain} with parameters specified as above and supplemented with the boundary conditions \eqref{Bd1}-\eqref{Bd2}.} of the value function $v(0,{\boldsymbol{\eta}})$ solving, at $t=0$, the HJB equation \eqref{hjb_ex} associated with \eqref{P_2D-proj}.
We observe that the value function  is smooth and convex, which justifies numerically the above assumption concerning the smoothness of the value function. 
Given the good control skills obtained from the controller synthesized from the reduced HJB equation  \eqref{hjb_ex} associated with \eqref{P_2D-proj}, one can reasonably infer that the corresponding ``reduced'' value function $v(0,\boldsymbol{\eta})$ shown in Fig.~\ref{fig:VF} provides thus a good approximation of the ``full'' value function associated with the DDE optimal control problem. Due to its simple shape, a polynomial approximation of $v(0,\boldsymbol{\eta})$ 
can be easily computed and gives (up to a small residual\footnote{The root mean square error between $v(0, \boldsymbol{\eta})$ and the polynomial approximation is $12\times 10^{-4}$.})
\be
v(0,\boldsymbol{\eta}) = 28.0025 \eta_1^2 - 7.8733 \eta_1 \eta_2 + 0.0046 \eta_1+ 10.6258 \eta_2^2 + 0.0525 \eta_2,
\ee
i.e.~a quadratic polynomial asymmetric in the $\eta_1$- and $\eta_2$-directions. 
Besides the quadratic cost functional used here, such a low-order approximation of the value function is probably due to the proximity to the first criticality and not limited to time $t=0$ shown here. This observation deserves of course more understanding which although not pursued here, is, we believe, worth mentioning.

\section{Concluding remarks} \label{Sect_conclusion}

As mentioned earlier, the good numerical skills shown in Sections \ref{Sec_num_results} and \ref{Sec_HJB}, as obtained from low-dimensional Galerkin systems, are conditioned to the distance from the first criticality, here the Hopf bifurcation.  In fact a numerical estimation of the RHS in the error estimate  \eqref{PM_value_est_goal} of Theorem~\ref{Thm_PM_val} reveals that the corresponding residual energies are actually small (not shown), in agreement with the numerical results.

As one gets further away from the first criticality, a larger number of Koornwinder polynomials is typically required to dispose of good GK approximations of, already,  the uncontrolled dynamics; see \cite{CGLW15}.  The numerical burden of the synthesis of controls at a nearly optimal cost---by solving the HJB equations corresponding to the relevant GK systems---becomes then quickly prohibitive, especially for the case of locally distributed controls\footnote{See \cite[Sect.~7]{CL15} for such a situation in the context of a viscous Burgers equation with instability.}.  One avenue to work within reduced state spaces of further reduced dimension compared to what would be required by a GK approximation, is to search for high-mode parameterizations that help reduce the residual energy contained in the unresolved modes, i.e.~to reduce a quantity like $ \mathcal{R}(u^\ast_{t,x},u^{N,\ast}_{t,x_{N}})$ in \eqref{PM_value_est_goal} that involves the residual energies, $\|\Pi_N^\perp y_{t,x}(\cdot; u^\ast_{t,x})\|_{L^2(t,T; D(\mathcal{A}))}$ and 
$\| \Pi_N^\perp y_{t,x}(\cdot; u^{N,\ast}_{t,x_{N}})\|_{L^2(t,T; D(\mathcal{A}))}$.

The theory of parameterizing manifolds (PM) \cite{CL15,CL16post,CLMcW2016} allows for such a reduction leading typically to approximate controls  coming essentially with error estimates that introduce a multiplying factor $0\leq Q <1$ (related to the PM) 
in an RHS similar to that of \eqref{PM_value_est_goal}; see \cite[Theorem 1 \& Corollary 2]{CL15}. The combination of the GK framework of \cite{CGLW15} with the PM reduction techniques of \cite{CL15} constitutes thus an idea that is worth pursuing for solving efficiently optimal control problems of nonlinear DDEs. 

\section*{Acknowledgments} 
This work has been partially supported by the Office of Naval Research (ONR) Multidisciplinary University Research Initiative (MURI) grant N00014-12-1-0911 and N00014-16-1-2073 (MDC), and by the National Science Foundation grants DMS-1616981 (MDC) and DMS-1616450 (HL). MDC and AK were supported by the project ``Optimal control of partial differential equations using parameterizing manifolds, model reduction, and dynamic programming'' funded by the {\it Foundation Hadamard/Gaspard Monge Program for Optimization and Operations Research (PGMO).}

\bibliographystyle{amsalpha}
\providecommand{\bysame}{\leavevmode\hbox to3em{\hrulefill}\thinspace}
\providecommand{\MR}{\relax\ifhmode\unskip\space\fi MR }
\providecommand{\MRhref}[2]{%
  \href{http://www.ams.org/mathscinet-getitem?mr=#1}{#2}
}
\providecommand{\href}[2]{#2}

\end{document}